\documentclass[a4paper,twoside,10pt]{amsart}
\usepackage{amsmath,amsfonts,amssymb,amsthm}
\newtheorem{theorem}{Theorem}[section]

\usepackage{mathrsfs}
\usepackage{color,graphicx}
\usepackage[a4paper]{geometry}
\numberwithin{equation}{section}

\author[G. Nemes]{Gerg\H{o} Nemes}
\address{Central European University, Department of Mathematics and its Applications, H-1051 Budapest, N\'ador utca 9, Hungary}
\email{nemesgery@gmail.com}


\begin{document}

\title[Resurgence of the incomplete gamma function]{The resurgence properties of\\ the incomplete gamma function II}

\begin{abstract} In this paper we derive a new representation for the incomplete gamma function, exploiting the reformulation of the method of steepest descents by C. J. Howls (Howls, Proc. R. Soc. Lond. A \textbf{439} (1992) 373--396). Using this representation, we obtain numerically computable bounds for the remainder term of the asymptotic expansion of the incomplete gamma function $\Gamma \left( { - a,\lambda a} \right)$ with large $a$ and fixed positive $\lambda$, and an asymptotic expansion for its late coefficients. We also give a rigorous proof of Dingle's formal result regarding the exponentially improved version of the asymptotic series of $\Gamma \left( { - a,\lambda a} \right)$.
\end{abstract}
\maketitle

\section{Introduction and main results}

It is known \cite[8.11.iii]{NIST}\cite{Temme} that, as $a\to \infty$ in the sector $\left| {\arg a} \right| \le  - \omega  + \frac{{3\pi }}{2} - \delta  <  - \omega  + \frac{{3\pi }}{2}$, for any fixed $\delta>0$, the incomplete gamma function has the following asymptotic expansion
\begin{equation}\label{eq10}
\Gamma \left( { - a,z} \right) \sim z^{ - a} e^{ - z} \sum\limits_{n = 0}^\infty  {\frac{{a^n b_n \left( { - \lambda } \right)}}{{\left( {z + a} \right)^{2n + 1} }}} ,
\end{equation}
with $z=\lambda a$, $\lambda>0$ and
\[
\omega := \omega\left(\lambda\right)=\arg \left( {\lambda  + \log \lambda  + 1 + \pi i} \right).
\]
The coefficients $b_n \left( \lambda  \right)$ originally appear in the related asymptotic series for $\Gamma \left( a,z \right)$ with $z = \lambda a$, $\lambda \geq 1$ \cite[8.11.E7]{NIST}. They are polynomials in $\lambda$ of degree $n$, the first few being
\[
b_0\left( \lambda  \right)=1,\; b_1\left( \lambda  \right)=\lambda, \; b_2\left( \lambda  \right)=2\lambda^2+\lambda,\; b_3 \left( \lambda  \right) = 6\lambda ^3  + 8\lambda ^2+\lambda.
\]
In the special case $a>0$, Dunster \cite{Dunster} and Gautschi \cite{Gautschi} gave computable error bounds for the series \eqref{eq10}. Dunster's error bound allows $\lambda$ to be a complex number with $\left|\arg \lambda\right|<\pi$. An asymptotic expansion for the coefficients $b_n \left( -\lambda  \right)$ was established by Dingle \cite[p. 162]{Dingle}. He also gave a re-expansion of the remainder term leading to an exponentially improved version of the asymptotic series \eqref{eq10} \cite[p. 463]{Dingle}. Nevertheless, the derivation of his results is based on formal interpretive, rather than rigorous, methods.

In the first part of this series of papers \cite{Nemes2}, we proved new resurgence-type representations for the remainder term of the asymptotic expansion of the incomplete gamma function $\Gamma \left( a,z \right)$ with $z = \lambda a$, $\lambda \geq 1$. Here resurgence has to be understood in the sense of Berry and Howls \cite{Berry}, meaning that the function (or a closely related function) reappears in the remainder of its own asymptotic series. These resurgence formulas for $\Gamma \left( a,z \right)$ have different forms according to whether $\lambda >1$ or $\lambda=1$. The main goal of this paper is to derive a similar representation for the remainder of the asymptotic expansion \eqref{eq10}. Our derivation is based on the reformulation of the method of steepest descents by Howls \cite{Howls}. Using this representation, we obtain numerically computable bounds for the remainder of the asymptotic series \eqref{eq10}, and an asymptotic expansion for its late coefficients. Our analysis also provides a rigorous treatment of Dingle's formal results.

In our investigations, the (phase) function $\omega$ plays a central role. It is a strictly decreasing function of $\lambda>0$. As $\lambda$ increases from $0$ to $+\infty$, $\omega$ decreases monotonically from $\pi$ to $0$.

As in the resurgence formulas found in the previous paper \cite{Nemes2}, the function that makes its appearance in the remainder term, is the scaled gamma function. This function is defined in terms of the classical gamma function as
\[
\Gamma^\ast \left( z \right) = \frac{{\Gamma \left( z \right)}}{{\sqrt {2\pi } z^{z - \frac{1}{2}} e^{ - z} }},
\]
for $\left|\arg z\right|<\pi$.

Our first theorem describes the resurgence properties of the incomplete gamma function $\Gamma \left( -a,z \right)$. The notations follow the ones given in \cite[8.11.iii]{NIST}.

\begin{theorem}\label{thm1} Let $\lambda >0$ be a fixed real number such that $z=\lambda a$, and let $N$ be a positive integer. Then
\begin{equation}\label{eq36}
\Gamma \left( { - a,z} \right) = z^{ - a} e^{ - z} \left( {\sum\limits_{n = 0}^{N - 1} {\frac{{a^n b_n \left( { - \lambda } \right)}}{{\left( {z + a} \right)^{2n + 1} }}}  + R_N \left( {a,\lambda } \right)} \right)
\end{equation}
for $\left| {\arg a} \right| <  - \omega  + \pi$, with
\begin{gather}\label{eq9}
\begin{split}
& b_n \left( { - \lambda } \right) = \left( {\lambda  + 1} \right)^n \left[ {\frac{{d^n }}{{dt^n }}\left( {\frac{{\left( {\lambda  + 1} \right)t}}{{\lambda e^t  + t - \lambda }}} \right)^{n + 1} } \right]_{t=0} \\ & = \left( { - 1} \right)^n \frac{{\left( {\lambda  + 1} \right)^{2n + 1} }}{{\sqrt {2\pi } i}}\int_0^{ + \infty } {t^{n - \frac{1}{2}} e^{ - t\left| {\lambda  + \log \lambda  + 1 + \pi i} \right|} \left( {e^{\left( {n + \frac{1}{2}} \right)\omega i} \Gamma^\ast \left( {te^{i\omega } } \right) - e^{ - \left( {n + \frac{1}{2}} \right)\omega i} \Gamma^\ast \left( {te^{ - i\omega } } \right)} \right)dt} ,
\end{split}
\end{gather}
where the second representation is true only for $n\geq 1$. The remainder term $R_N \left( a,\lambda \right)$ is given by
\begin{gather}\label{eq8}
\begin{split}
R_N \left( {a,\lambda } \right) = \; & \frac{{\left( { - 1} \right)^N a^N }}{{\left( {z + a} \right)^{2N + 1} }}\left( {\lambda  + 1} \right)^{2N + 1} \frac{{e^{ \left( {N + \frac{1}{2}} \right)\omega i} }}{{\sqrt {2\pi } i}}\int_0^{ + \infty } {\frac{{t^{N - \frac{1}{2}} e^{ - t\left| {\lambda  + \log \lambda  + 1 + \pi i} \right|} }}{{1 + te^{i\omega } /a}}\Gamma^\ast \left( {te^{i\omega } } \right)dt} \\
\\ & + \frac{{\left( { - 1} \right)^{N + 1} a^N }}{{\left( {z + a} \right)^{2N + 1} }}\left( {\lambda  + 1} \right)^{2N + 1} \frac{{e^{ - \left( {N + \frac{1}{2}} \right) \omega i } }}{{\sqrt {2\pi } i}}\int_0^{ + \infty } {\frac{{t^{N - \frac{1}{2}} e^{ - t\left| {\lambda  + \log \lambda  + 1 + \pi i} \right|} }}{{1 + te^{ - i\omega } /a}}\Gamma^\ast \left( {te^{ - i\omega } } \right)dt} .
\end{split}
\end{gather}
\end{theorem}

In Section \ref{section3}, we will show how to obtain numerically computable bounds for the remainder term $R_N \left( {a,\lambda } \right)$ using its explicit form given in Theorem \ref{thm1}. Some other formulas for the coefficients $b_n\left(-\lambda\right)$ can be found in \cite[Appendix A]{Nemes2}.

It is seen from the monotonicity properties of $\omega$ that the region of validity of our resurgence formula becomes wider as $\lambda$ becomes larger.

In deriving the further results, which are based on Theorem \ref{thm1}, it is required to estimate $\Gamma^\ast \left( {te^{\pm i\omega } } \right)$ or the remainder term in its asymptotic series (see below). However, as $\lambda$ approaches $0$ or, equivalently, $\omega$ approaches $\pi$, the simple estimates for these functions break down due to the presence of the poles of the gamma function along the negative real axis. Therefore, our results are less effective when $\lambda$ is close to $0$. While proving Theorem \ref{thm1}, we have found the following alternative representation for the remainder $R_N \left( {a,\lambda } \right)$ in the range $0<\lambda <W\left(e^{-1}\right)=0.27846\ldots$, or equivalently $\frac{\pi}{2} < \omega <\pi$. Here $W$ denotes the principal branch of the Lambert $W$-function \cite[4.13]{NIST}. In this representation, the scaled gamma function is evaluated in the right-half plane, in the region where we have well-behaved simple estimates for it. Thus, it is possible to derive alternative estimations based on this representation at the cost of having infinite series in the final results. We do not discuss the details in the present paper.

\begin{theorem}\label{thm4} Let $0<\lambda <W\left(e^{-1}\right)=0.27846\ldots$ be a fixed real number, and let $N$ be a positive integer. For any non-negative integer $k$ let $\omega _k  = \arg \left( {\lambda  + \log \lambda  + 1 + \left( {2k + 1} \right)\pi i} \right)$ (note that $\omega_0=\omega$). Then the remainder term $R_N \left( {a,\lambda } \right)$ defined by \eqref{eq36} has the expansion
\begin{align*}
R_N \left( {a,\lambda } \right) = \; & \frac{{\left( { - 1} \right)^N a^N }}{{\left( {z + a} \right)^{2N + 1} }}\left( {\lambda  + 1} \right)^{2N + 1} \sum\limits_{k = 0}^\infty  {\frac{{e^{\left( {N + \frac{1}{2}} \right)\omega _k i} }}{{\sqrt {2\pi } i}}\int_0^{ + \infty } {\frac{{t^{N - \frac{1}{2}} e^{ - t\left| {\lambda  + \log \lambda  + 1 - \left( {2k + 1} \right)\pi i} \right|} }}{{1 + te^{i\omega _k } /a}}\frac{{dt}}{{\Gamma^\ast \left( {te^{i\left( {\omega _k  - \pi } \right)} } \right)}}} } 
\\ & + \frac{{\left( { - 1} \right)^{N + 1} a^N }}{{\left( {z + a} \right)^{2N + 1} }}\left( {\lambda  + 1} \right)^{2N + 1} \sum\limits_{k = 0}^\infty  {\frac{{e^{ - \left( {N + \frac{1}{2}} \right)\omega _k i} }}{{\sqrt {2\pi } i}}\int_0^{ + \infty } {\frac{{t^{N - \frac{1}{2}} e^{ - t\left| {\lambda  + \log \lambda  + 1 + \left( {2k + 1} \right)\pi i} \right|} }}{{1 + te^{ - i\omega _k } /a}}\frac{{dt}}{{\Gamma^\ast \left( {te^{ - i\left( {\omega _k  - \pi } \right)} } \right)}}} },
\end{align*}
for $\left| {\arg a} \right| <  - \omega  + \pi$.
\end{theorem}

It is interesting to note that this representation for $R_N \left( {a,\lambda } \right)$ has infinitely many different singular directions giving rise to infinitely many Stokes lines, whereas the representation \eqref{eq8} produces only two Stokes lines. This phenomenon, namely when the number of Stokes lines depends on a certain parameter (in our case on $\lambda$) is related to the higher-order Stokes phenomenon (see, e.g., \cite{Howls2}). However, we will see in Section \ref{section3} that, in our case, the infinitely many contributions can be summed up explicitly leading to the expression \eqref{eq8}, ending up with only two Stokes lines. The critical value $\lambda = W\left(e^{-1}\right)$ (corresponding to $\omega = \frac{\pi}{2}$) and the infinite sums in Theorem \ref{thm4} are consequences of the fact that the rays $\left|\arg z\right| = \frac{\pi}{2}$ are the Stokes lines for the function $\Gamma^\ast \left( z\right)$, and that along these rays infinitely many exponentially small terms appear in the asymptotic expansion of $\Gamma^\ast \left( z\right)$ (see, e.g., \cite{Nemes}).

In the important paper \cite{Boyd}, Boyd gave the following resurgence formula for the well-known asymptotic expansion of scaled gamma function (see also \cite{Nemes}):
\begin{equation}\label{eq11}
\Gamma^\ast \left( z \right) = \sum\limits_{n = 0}^{N - 1} {\left(-1\right)^n\frac{\gamma_n }{z^n }}  + M_N \left( z \right),
\end{equation}
for $\left|\arg z\right| < \frac{\pi}{2}$ and $N\geq 1$, with
\begin{equation}\label{eq14}
M_N \left( z \right) = \frac{1}{{2\pi i}}\frac{i^N}{{z^N }}\int_0^{ + \infty } {\frac{{t^{N - 1} e^{ - 2\pi t} \Gamma^\ast \left( {it} \right)}}{{1 - it/z}}dt}  - \frac{1}{{2\pi i}}\frac{{\left( { - i} \right)^N }}{{z^N }}\int_0^{ + \infty } {\frac{{t^{N - 1} e^{ - 2\pi t} \Gamma^\ast \left( { - it} \right)}}{{1 + it/z}}dt} .
\end{equation}
Here the $\gamma_n$'s are the so-called Stirling coefficients (see \cite{Nemes}). This representation of the scaled gamma function will play an important role in later sections of this paper. By \eqref{eq11}, $M_N \left( z \right)$ can be defined in the wider range $\left|\arg z\right| <\pi$, and it is known that $M_N \left( z \right)=\mathcal{O}\left(\left|z\right|^{-N}\right)$ as $z \to \infty$ in the sector $\left|\arg z\right|\leq \pi-\delta< \pi$, for any fixed $\delta >0$.

The rest of the paper is organised as follows. In Section \ref{section2}, we prove the resurgence formulas stated in Theorems \ref{thm1} and \ref{thm4}. In Section \ref{section3}, we give explicit and numerically computable error bounds for the asymptotic series \eqref{eq10} using the results of Theorem \ref{thm1}. In Section \ref{section4}, asymptotic approximations for the coefficients $b_n \left( -\lambda  \right)$ are given. In Section \ref{section5}, we formulate and prove a rigorous form of Dingle's exponentially improved version of the asymptotic expansion \eqref{eq10}. 

\section{Proof of the resurgence formula}\label{section2} 
\begin{figure}[!t]
\def\svgwidth{0.6\columnwidth}
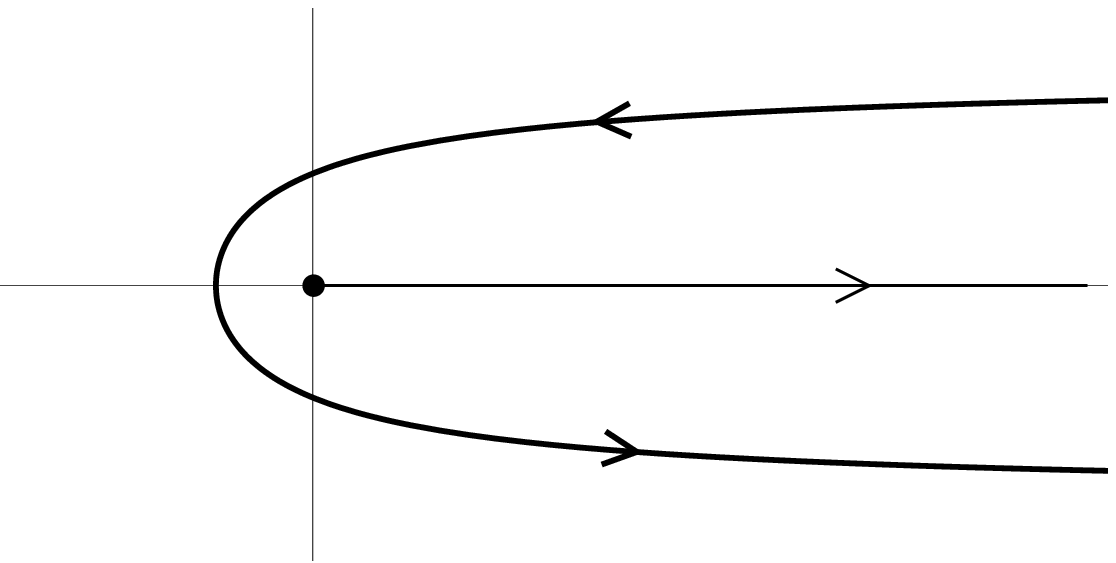
\caption{The contour $\Gamma$ encircling the path $\mathscr{P}\left(0\right)$.}
\label{fig1}
\end{figure}

Our analysis is based on the following integral representation (see, e.g., \cite[8.6.E7]{NIST})
\[
\Gamma \left( { - a,z} \right) = z^{ - a} \int_0^{ + \infty } {\exp \left( { - at - ze^t } \right)dt} \; \text{ for } \; \Re\left(z\right)>0.
\]
If $z=\lambda a$, with $\lambda>0$ fixed, then
\begin{equation}\label{eq1}
\Gamma \left( { - a,z} \right) = z^{ - a} e^{ - z} \int_0^{ + \infty } {\exp \left( { - a\left( {\lambda e^t  + t - \lambda } \right)} \right)dt} 
\end{equation}
provided that $\Re\left(a\right)>0$. The saddle points of the integrand are the roots of the equation $\lambda e^t  + 1=0$. Hence, the saddle points are given by $t^{\left(k\right)}=-\log \lambda +\left(2k+1\right)\pi i$ where $k$ is an arbitrary integer. We denote by $\mathscr{C}^{\left(k\right)}\left(\theta\right)$ the portion of the steepest paths that pass through the saddle point $t^{\left(k\right)}$. Here, and subsequently, we write $\theta = \arg a$. As for the path of integration $\mathscr{P}\left(\theta\right)$ in \eqref{eq1}, we shall take that connected component of
\[
\left\{ {t \in \mathbb{C}:\arg \left[ {e^{i\theta } \left( {\lambda e^t + t - \lambda} \right)} \right] = 0} \right\} ,
\]
which is the positive real axis for $\theta=0$ and is the continuous deformation of the positive real axis as $\theta$ varies. Let $f\left( {t,\lambda } \right) = \lambda e^t  + t - \lambda$ for some fixed $\lambda > 0$. Hence, we can write \eqref{eq1} as
\[
\Gamma \left( { - a,z} \right) = z^{ - a} e^{ - z} \int_0^{ + \infty } {e^{ - af\left( {t,\lambda } \right)} dt} 
\]
for any $\Re \left( a \right) > 0$. For simplicity, we assume that $a>0$. In due course, we shall appeal to an analytic continuation argument to extend our results to complex $a$. If
\begin{equation}\label{eq2}
\tau = f\left( {t,\lambda } \right),
\end{equation}
then $\tau$ is real on the curve $\mathscr{P}\left(0\right)$, and, as $t$ travels along this curve from $0$ to $+\infty$, $\tau$ increases from $0$ to $+\infty$. Therefore, corresponding to each positive value of $\tau$, there is a value of $t$, say $t\left(\tau\right)$, satisfying \eqref{eq2} with $t\left(\tau\right)>0$. In terms of $\tau$, we have
\[
\Gamma \left( { - a,z} \right) = z^{ - a} e^{ - z} \int_0^{ + \infty } {e^{ - a\tau } \frac{{dt}}{{d\tau }}d\tau }  = z^{ - a} e^{ - z} \int_0^{ + \infty } {\frac{{e^{ - a\tau } }}{{f'\left( {t\left( \tau  \right),\lambda } \right)}}d\tau } .
\]
Following Howls \cite{Howls}, we express the function involving $t\left(\tau\right)$ as a contour integral using the residue theorem, to find
\[
\Gamma \left( { - a,z} \right) = z^{ - a} e^{ - z} \int_0^{ + \infty } {e^{ - a\tau } \frac{1}{{2\pi i}}\oint_\Gamma  {\frac{{f^{ - 1} \left( {u,\lambda } \right)}}{{1 - \tau f^{ - 1} \left( {u,\lambda } \right)}}du} d\tau } ,
\]
where the contour $\Gamma$ encircles the path $\mathscr{P}\left(0\right)$ in the positive direction and does not enclose any of the saddle points $t^{\left(k\right)}$ (see Figure \ref{fig1}). Now, we employ the well-known expression for non-negative integer $N$
\[
\frac{1}{1 - x} = \sum\limits_{n = 0}^{N-1} {x^n}  + \frac{x^N}{1 - x},\; x \neq 1,
\]
to expand the function under the contour integral in powers of $\tau f^{ - 1} \left( {u,\lambda } \right)$. The result is
\[
\Gamma \left( { - a,z} \right) = z^{ - a} e^{ - z} \left( {\sum\limits_{n = 0}^{N - 1} {\int_0^{ + \infty } {\tau ^n e^{ - a\tau } \frac{1}{{2\pi i}}\oint_\Gamma  {\frac{{du}}{{f^{n + 1} \left( {u,\lambda } \right)}}} d\tau } }  + R_N \left( {a,\lambda } \right)} \right)
\]
where
\begin{equation}\label{eq3}
R_N \left( {a,\lambda } \right) = \int_0^{ + \infty } {\tau ^N e^{ - a\tau } \frac{1}{{2\pi i}}\oint_\Gamma  {\frac{{f^{ - N - 1} \left( {u,\lambda } \right)}}{{1 - \tau f^{ - 1} \left( {u,\lambda } \right)}}du} d\tau } .
\end{equation}
The path $\Gamma$ in the sum can be shrunk into a small circle around $0$, and we arrive at
\begin{equation}\label{eq4}
\Gamma \left( { - a,z} \right) = z^{ - a} e^{ - z} \left( {\sum\limits_{n = 0}^{N - 1} {\frac{{a^n b_n \left( { - \lambda } \right)}}{{\left( {z + a} \right)^{2n + 1} }}}  + R_N \left( {a,\lambda } \right)} \right),
\end{equation}
where
\[
b_n \left( { - \lambda } \right) = \frac{{\left( {\lambda  + 1} \right)^{2n + 1} \Gamma \left( {n + 1} \right)}}{{2\pi i}}\oint_{\left(0^+\right)}  {\frac{{du}}{{f^{n + 1} \left( {u,\lambda } \right)}}}  = \left( {\lambda  + 1} \right)^n \left[ {\frac{{d^n }}{{dt^n }}\left( {\frac{{\left( {\lambda  + 1} \right)t}}{{\lambda e^t  + t - \lambda }}} \right)^{n + 1} } \right]_{t=0} .
\]

\begin{figure}[!t]
\def\svgwidth{0.5\columnwidth}
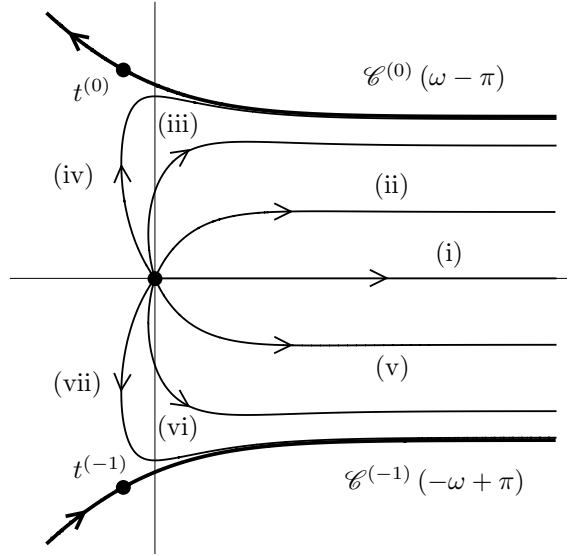
\caption{The path $\mathscr{P}\left(\theta\right)$ emanating from the origin when $\lambda=2$ and (i) $\theta=0$, (ii) $\theta=-1$, (iii) $\theta=-2$, (iv) $\theta=-\frac{12}{5}$, (v) $\theta=1$, (vi) $\theta=2$ and (vii) $\theta=\frac{12}{5}$. The paths $\mathscr{C} ^{\left( 0 \right)} \left( \omega-\pi \right)$ and $\mathscr{C}^{\left( -1 \right)} \left( -\omega+\pi \right)$ are the adjacent contours to $0$.}
\label{fig2}
\end{figure}

Performing the change of variable $a\tau  = s$ in \eqref{eq3} yields
\begin{equation}\label{eq5}
R_N \left( {a,\lambda } \right) = \frac{{a^N }}{{\left( {z + a} \right)^{2N + 1} }}\left( {\lambda  + 1} \right)^{2N + 1} \int_0^{ + \infty } {s^N e^{ - s} \frac{1}{{2\pi i}}\oint_\Gamma  {\frac{{f^{ - N - 1} \left( {u,\lambda } \right)}}{{1 - \left( {s/a} \right)f^{ - 1} \left( {u,\lambda } \right)}}du} ds} 
\end{equation}
This representation of $R_N \left( {a,\lambda } \right)$ and the formula \eqref{eq4} can be continued analytically if we choose $\Gamma = \Gamma\left(\theta\right)$ to be an infinite contour that surrounds the path $\mathscr{P}\left(\theta\right)$ in the anti-clockwise direction and that does not encircle any of the saddle points $t^{\left(k\right)}$. This continuation argument works until the path $\mathscr{P}\left(\theta\right)$ runs into a saddle point. In the terminology of Howls, such saddle points are called adjacent to the endpoint $0$. At this point the analysis splits into two parts, because the number of adjacent saddles depends on the size of $\lambda$ (or $\omega$). First assume that $\lambda \geq W\left( {e^{ - 1} } \right) =0.27846\ldots$, where $W$ denotes the principal branch of the Lambert $W$-function \cite[4.13]{NIST}. Under this assumption, we have $0 < \omega \leq \frac{\pi }{2}$. In this case, when $\theta =  - \arg \left( { - \lambda  - \log \lambda  - 1 + \pi i} \right) = \arg \left( {\lambda  + \log \lambda  + 1 + \pi i} \right) - \pi =\omega -\pi$, the path $\mathscr{P}\left(\theta\right)$ connects to the saddle point $t^{\left(0\right)}=-\log \lambda+\pi i$; and when $\theta =  - \arg \left( { - \lambda  - \log \lambda  - 1 - \pi i} \right) = \arg \left( {\lambda  + \log \lambda  + 1 - \pi i} \right) + \pi = -\omega+\pi$, the path $\mathscr{P}\left(\theta\right)$ connects to the saddle point $t^{\left(-1\right)}=-\log \lambda-\pi i$. These are the adjacent saddles. The set
\[
\Delta = \left\{u\in \mathscr{P}\left(\theta\right):  \omega-\pi < \theta  <  -\omega+\pi \right\}
\]
forms a domain in the complex plane whose boundary contains portions of steepest descent paths through the adjacent saddles (see Figure \ref{fig2}). These paths are $\mathscr{C} ^{\left( 0 \right)} \left( \omega-\pi \right)$ and $\mathscr{C}^{\left( -1 \right)} \left( -\omega+\pi \right)$, and they are called the adjacent contours to the endpoint $0$. If $N\geq 1$, the function under the second integral sign in \eqref{eq5} is an analytic function of $u$ in the domain $\Delta$ and at the points between the adjacent contours, therefore we can deform $\Gamma$ over the adjacent contours. We thus find that for $\omega-\pi < \theta < -\omega+\pi$ and $N\geq 1$, \eqref{eq5} may be written
\begin{gather}\label{eq6}
\begin{split}
R_N \left( {a,\lambda } \right) = \; & \frac{{a^N }}{{\left( {z + a} \right)^{2N + 1} }}\left( {\lambda  + 1} \right)^{2N + 1} \int_0^{ + \infty } {s^N e^{ - s} \frac{1}{{2\pi i}}\int_{\mathscr{C} ^{\left( 0 \right)} \left( \omega-\pi \right)}  {\frac{{f^{ - N - 1} \left( {u,\lambda } \right)}}{{1 - \left( {s/a} \right)f^{ - 1} \left( {u,\lambda } \right)}}du} ds} \\ & + \frac{{a^N }}{{\left( {z + a} \right)^{2N + 1} }}\left( {\lambda  + 1} \right)^{2N + 1} \int_0^{ + \infty } {s^N e^{ - s} \frac{1}{{2\pi i}}\int_{\mathscr{C} ^{\left( -1 \right)} \left( -\omega+\pi \right)}  {\frac{{f^{ - N - 1} \left( {u,\lambda } \right)}}{{1 - \left( {s/a} \right)f^{ - 1} \left( {u,\lambda } \right)}}du} ds}.
\end{split}
\end{gather}
Now we make the change of variable
\[
s = t\frac{{\left| {f\left( { - \log \lambda  + \pi i,\lambda } \right) - f\left( {0,\lambda } \right)} \right|}}{{f\left( { - \log \lambda  + \pi i,\lambda } \right) - f\left( {0,\lambda } \right)}}f\left( {u,\lambda } \right) =  - te^{i\omega } f\left( {u,\lambda } \right)
\]
in the first, and
\[
s = t\frac{{\left| {f\left( { - \log \lambda  - \pi i,\lambda } \right) - f\left( {0,\lambda } \right)} \right|}}{{f\left( { - \log \lambda  - \pi i,\lambda } \right) - f\left( {0,\lambda } \right)}}f\left( {u,\lambda } \right) =  - te^{ - i\omega } f\left( {u,\lambda } \right)
\]
in the second double integral. Clearly, by the definition of the adjacent contours, $t$ is positive. The quantities $f\left( { - \log \lambda  \pm \pi i,\lambda } \right) - f\left( {0,\lambda } \right) =  - \lambda  - \log \lambda  - 1 \pm \pi i$ were essentially called ``singulants" by Dingle \cite[p. 147]{Dingle}. With this change of variable, the representation \eqref{eq6} for $R_N \left( {a,\lambda } \right)$ becomes
\begin{gather}\label{eq7}
\begin{split}
R_N \left( {a,\lambda } \right) = \; & \frac{{a^N }}{{\left( {z + a} \right)^{2N + 1} }}\left( {\lambda  + 1} \right)^{2N + 1} \left( { - e^{i\omega } } \right)^{N + 1} \int_0^{ + \infty } {\frac{{t^N }}{{1 + te^{i\omega } /a}}\frac{1}{{2\pi i}}\int_{\mathscr{C}^{\left( 0 \right)} \left( \omega-\pi  \right)} {e^{te^{i\omega } f\left( {u,\lambda } \right)} du} dt} 
\\& +\frac{{a^N }}{{\left( {z + a} \right)^{2N + 1} }}\left( {\lambda  + 1} \right)^{2N + 1} \left( { - e^{ - i\omega } } \right)^{N + 1} \int_0^{ + \infty } {\frac{{t^N }}{{1 + te^{ - i\omega } /a}}\frac{1}{{2\pi i}}\int_{\mathscr{C}^{\left( { - 1} \right)} \left( { - \omega+\pi } \right)} {e^{te^{ - i\omega } f\left( {u,\lambda } \right)} du} dt} ,
\end{split}
\end{gather}
for $\omega-\pi < \theta  <  -\omega+\pi$ and $N\geq 1$. It remains to evaluate the contour integrals. The change of variable $u = -\log\lambda + v + \pi i$ yields
\[
\frac{1}{{2\pi i}}\int_{\mathscr{C}^{\left( 0 \right)} \left( \omega -\pi \right)} {e^{te^{i\omega } f\left( {u,\lambda } \right)} du}  = \frac{{e^{ - t\left| {\lambda  + \log \lambda  + 1 - \pi i} \right|} }}{{2\pi i}}\int_{\mathscr{L}\left( \omega  \right)} {e^{ - te^{i\omega } \left( {e^v  - v - 1} \right)} dv} ,
\]
where $\mathscr{L}\left( \omega  \right)$ is the image of $\mathscr{C}^{\left( 0 \right)} \left( \omega -\pi \right)$ in the $v$-plane. Along $\mathscr{L}\left( \omega  \right)$, it holds that $\arg \left[ {e^{i\omega } \left( {e^v  - v - 1} \right)} \right] = 0$, whence $\mathscr{L}\left( \omega  \right)$ is the steepest descent path through the saddle point at $v=0$. Since $\lambda \geq  W\left( {e^{ - 1} } \right)$, i.e., $0<\omega \leq \frac{\pi}{2}$, the contour integral can be expressed in terms of the scaled gamma function (see \cite[equation (2.11)]{Boyd}) leading to the formula
\[
\frac{1}{{2\pi i}}\int_{\mathscr{C}^{\left( 0 \right)} \left( {\omega  - \pi } \right)} {e^{te^{i\omega } f\left( {u,\lambda } \right)} du}  =- \frac{{e^{ - t\left| {\lambda  + \log \lambda  + 1 - \pi i} \right|} }}{{\sqrt {2\pi } i}}t^{ - \frac{1}{2}} e^{ - \frac{\omega }{2}i} \Gamma^\ast \left( {te^{i\omega } } \right).
\]
Similarly, we find that
\[
\frac{1}{{2\pi i}}\int_{\mathscr{C}^{\left( { - 1} \right)} \left( {-\omega  + \pi } \right)} {e^{te^{ - i\omega } f\left( {u,\lambda } \right)} du}  = \frac{{e^{ - t\left| {\lambda  + \log \lambda  + 1 + \pi i} \right|} }}{{\sqrt {2\pi } i}}t^{ - \frac{1}{2}} e^{\frac{\omega }{2}i} \Gamma^\ast  \left( {te^{ - i\omega } } \right) .
\]
Substituting these expressions into \eqref{eq7} and using the fact that $\left| {\lambda  + \log \lambda  + 1 - \pi i} \right| = \left| {\lambda  + \log \lambda  + 1 + \pi i} \right|$, gives \eqref{eq8} for $\lambda \geq W\left( {e^{ - 1} } \right)$.

\begin{figure}[!t]
\def\svgwidth{0.55\columnwidth}
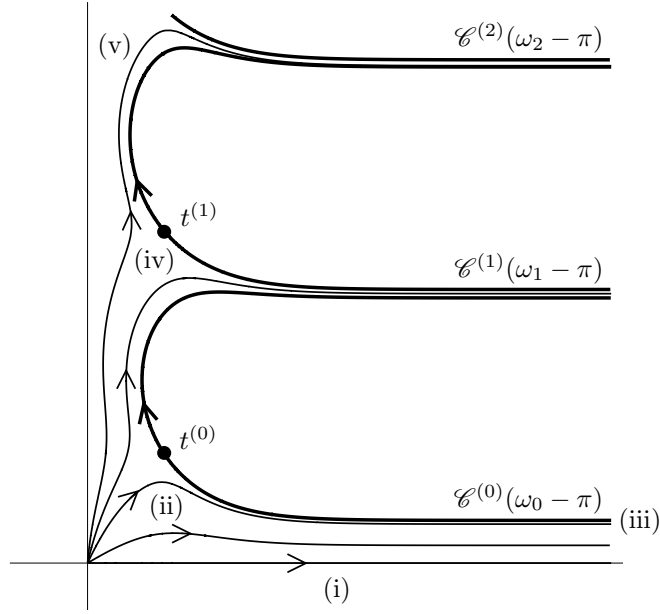
\caption{The path $\mathscr{P}\left(\theta\right)$ emanating from the origin when $\lambda=\frac{1}{10}$ and (i) $\theta=0$, (ii) $\theta=-\frac{1}{2}$, (iii) $\theta=-\frac{11}{10}$, (iv) $\theta=-\frac{53}{40}$ and (v) $\theta=-\frac{37}{25}$. The paths $\mathscr{C} ^{\left( k \right)} \left( \omega_k-\pi \right)$ are adjacent contours to $0$.}
\label{fig3}
\end{figure}

Consider now the remaining interval $0<\lambda< W\left( {e^{ - 1} } \right)$. In this case, all the infinitely many saddles are adjacent to the endpoint $0$. For any non-negative integer $k$, let $\omega _k  = \arg \left( {\lambda  + \log \lambda  + 1 + \left( {2k + 1} \right)\pi i} \right)$. Note that $\omega_0=\omega$. When $\theta =\omega_k -\pi$, the path $\mathscr{P}\left(\theta\right)$ connects to the saddle point $t^{\left(k\right)}=-\log \lambda+\left(2k+1\right)\pi i$; and when $\theta = -\omega_k+\pi$, the path $\mathscr{P}\left(\theta\right)$ connects to the saddle point $t^{\left(-k-1\right)}=-\log \lambda-\left(2k+1\right)\pi i$. The set
\[
\Delta = \left\{u\in \mathscr{P}\left(\theta\right):  \omega-\pi < \theta  <  -\omega+\pi,\; \omega _{k + 1}  - \pi  <  \pm \theta  < \omega _k  - \pi,\; k\geq 0 \right\}
\]
forms a domain in the complex plane whose boundary contains steepest descent paths through the adjacent saddles (see Figure \ref{fig3}). These paths are $\mathscr{C} ^{\left( k \right)} \left( \omega_k-\pi \right)$ and $\mathscr{C}^{\left( -k-1 \right)} \left( -\omega_k+\pi \right)$, the adjacent contours to the endpoint $0$. If $N\geq 1$, the function under the second integral sign in \eqref{eq5} is an analytic function of $u$ in the domain $\Delta$ and at the points between the adjacent contours, therefore we can deform $\Gamma$ over the adjacent contours. We thus find that for $\omega-\pi < \theta < -\omega+\pi$ and $N\geq 1$, \eqref{eq5} may be written
\begin{gather}\label{eq30}
\begin{split}
R_N \left( {a,\lambda } \right) = \; & \frac{{a^N }}{{\left( {z + a} \right)^{2N + 1} }}\left( {\lambda  + 1} \right)^{2N + 1} \int_0^{ + \infty } {s^N e^{ - s} \sum\limits_{k = 0}^\infty  {\frac{1}{{2\pi i}}\int_{\mathscr{C}^{\left( k \right)} \left( {\omega _k  - \pi } \right)} {\frac{{f^{ - N - 1} \left( {u,\lambda } \right)}}{{1 - \left( {s/a} \right)f^{ - 1} \left( {u,\lambda } \right)}}du} } ds} \\
& + \frac{{a^N }}{{\left( {z + a} \right)^{2N + 1} }}\left( {\lambda  + 1} \right)^{2N + 1} \int_0^{ + \infty } {s^N e^{ - s} \sum\limits_{k = 0}^\infty  {\frac{1}{{2\pi i}}\int_{\mathscr{C}^{\left( { - k - 1} \right)} \left( { - \omega _k  + \pi } \right)} {\frac{{f^{ - N - 1} \left( {u,\lambda } \right)}}{{1 - \left( {s/a} \right)f^{ - 1} \left( {u,\lambda } \right)}}du} } ds} .
\end{split}
\end{gather}

Now we make the changes of variable
\[
s = t\frac{{\left| {f\left( { - \log \lambda  + \left( {2k + 1} \right)\pi i,\lambda } \right) - f\left( {0,\lambda } \right)} \right|}}{{f\left( { - \log \lambda  + \left( {2k + 1} \right)\pi i,\lambda } \right) - f\left( {0,\lambda } \right)}}f\left( {u,\lambda } \right) =  - te^{i\omega _k } f\left( {u,\lambda } \right)
\]
in the integrals under the first sum, and
\[
s = t\frac{{\left| {f\left( { - \log \lambda  - \left( {2k + 1} \right)\pi i,\lambda } \right) - f\left( {0,\lambda } \right)} \right|}}{{f\left( { - \log \lambda  - \left( {2k + 1} \right)\pi i,\lambda } \right) - f\left( {0,\lambda } \right)}}f\left( {u,\lambda } \right) =  - te^{ - i\omega _k } f\left( {u,\lambda } \right)
\]
in the integrals under the second sum. By the definition of the adjacent contours, $t$ is positive. In this case, Dingle's singulants are $f\left( { - \log \lambda  + \left( {2k + 1} \right)\pi i,\lambda } \right) - f\left( {0,\lambda } \right) =  - \lambda  - \log \lambda  - 1 + \left( {2k + 1} \right)\pi i$ and $f\left( { - \log \lambda  - \left( {2k + 1} \right)\pi i,\lambda } \right) - f\left( {0,\lambda } \right) =  - \lambda  - \log \lambda  - 1 - \left( {2k + 1} \right)\pi i$. With these changes of variable, the representation \eqref{eq30} for $R_N \left( {a,\lambda } \right)$ becomes
\begin{gather}\label{eq31}
\begin{split}
& R_N \left( {a,\lambda } \right) = \frac{{a^N }}{{\left( {z + a} \right)^{2N + 1} }}\left( {\lambda  + 1} \right)^{2N + 1} \sum\limits_{k = 0}^\infty  {\left( { - e^{i\omega _k } } \right)^{N + 1} \int_0^{ + \infty } {\frac{{t^N }}{{1 + te^{i\omega _k } /a}}\frac{1}{{2\pi i}}\int_{\mathscr{C}^{\left( k \right)} \left( {\omega _k  - \pi } \right)} {e^{te^{i\omega _k } f\left( {u,\lambda } \right)} du} dt} } 
\\ & + \frac{{a^N }}{{\left( {z + a} \right)^{2N + 1} }}\left( {\lambda  + 1} \right)^{2N + 1} \sum\limits_{k = 0}^\infty  {\left( { - e^{ - i\omega _k } } \right)^{N + 1} \int_0^{ + \infty } {\frac{{t^N }}{{1 + te^{ - i\omega _k } /a}}\frac{1}{{2\pi i}}\int_{\mathscr{C}^{\left( { - k - 1} \right)} \left( { - \omega _k  + \pi } \right)} {e^{te^{ - i\omega _k } f\left( {u,\lambda } \right)} du} dt} }, 
\end{split}
\end{gather}
for $\omega -\pi <\theta<-\omega+\pi$ and $N \geq 1$. Now we evaluate the contour integrals. The change of variable $u =  - \log \lambda  + v + \left( {2k + 1} \right)\pi i$ gives
\[
\frac{1}{{2\pi i}}\int_{\mathscr{C}^{\left( k \right)} \left( {\omega _k  - \pi } \right)} {e^{te^{i\omega _k } f\left( {u,\lambda } \right)} du}  = \frac{{e^{ - t\left| {\lambda  + \log \lambda  + 1 - \left( {2k + 1} \right)\pi i} \right|} }}{{2\pi i}}\int_{\mathscr{L}\left( {\omega _k } \right)} {e^{te^{i\left( {\omega _k  - \pi } \right)} \left( {e^v  - v - 1} \right)} dv} ,
\]
where $\mathscr{L}\left( {\omega _k } \right)$ is the image of $\mathscr{C}^{\left( k \right)} \left( {\omega _k  - \pi } \right)$ in the $v$-plane. Along $\mathscr{L}\left( {\omega _k } \right)$, it holds that $\arg \left[ { - e^{i\left( {\omega _k  - \pi } \right)} \left( {e^v  - v - 1} \right)} \right] = 0$, whence $\mathscr{L}\left( {\omega _k } \right)$ is the steepest descent path through the saddle point at $v = 0$. Since  $0<\lambda< W\left( {e^{ - 1} } \right)$, i.e., $\frac{\pi}{2}<\omega_k\leq \omega_0=\omega <\pi$, the contour integrals can be expressed in terms of the reciprocal of the scaled gamma function leading to the formula
\[
\frac{1}{{2\pi i}}\int_{\mathscr{C}^{\left( k \right)} \left( {\omega _k  - \pi } \right)} {e^{te^{i\omega _k } f\left( {u,\lambda } \right)} du}  =  - \frac{{e^{ - t\left| {\lambda  + \log \lambda  + 1 - \left( {2k + 1} \right)\pi i} \right|} }}{{\sqrt {2\pi } i}}\frac{{t^{ - \frac{1}{2}} e^{ - \frac{{\omega _k }}{2}i} }}{{\Gamma^\ast \left( {te^{i\left( {\omega _k  - \pi } \right)} } \right)}}.
\]
Similarly, we find that
\[
\frac{1}{{2\pi i}}\int_{\mathscr{C}^{\left( { - k - 1} \right)} \left( { - \omega _k  + \pi } \right)} {e^{te^{ - i\omega _k } f\left( {u,\lambda } \right)} du}  = \frac{{e^{ - t\left| {\lambda  + \log \lambda  + 1 + \left( {2k + 1} \right)\pi i} \right|} }}{{\sqrt {2\pi } i}}\frac{{t^{ - \frac{1}{2}} e^{\frac{{\omega _k }}{2}i} }}{{\Gamma^\ast \left( {te^{ - i\left( {\omega _k  - \pi } \right)} } \right)}}.
\]
Substituting these expressions into \eqref{eq31}, yields
\begin{align*}
R_N \left( {a,\lambda } \right) = \; & \frac{{\left( { - 1} \right)^N a^N }}{{\left( {z + a} \right)^{2N + 1} }}\left( {\lambda  + 1} \right)^{2N + 1} \sum\limits_{k = 0}^\infty  {\frac{{e^{\left( {N + \frac{1}{2}} \right)\omega _k i} }}{{\sqrt {2\pi } i}}\int_0^{ + \infty } {\frac{{t^{N - \frac{1}{2}} e^{ - t\left| {\lambda  + \log \lambda  + 1 - \left( {2k + 1} \right)\pi i} \right|} }}{{1 + te^{i\omega _k } /a}}\frac{{dt}}{{\Gamma^\ast \left( {te^{i\left( {\omega _k  - \pi } \right)} } \right)}}} } 
\\ & + \frac{{\left( { - 1} \right)^{N + 1} a^N }}{{\left( {z + a} \right)^{2N + 1} }}\left( {\lambda  + 1} \right)^{2N + 1} \sum\limits_{k = 0}^\infty  {\frac{{e^{ - \left( {N + \frac{1}{2}} \right)\omega _k i} }}{{\sqrt {2\pi } i}}\int_0^{ + \infty } {\frac{{t^{N - \frac{1}{2}} e^{ - t\left| {\lambda  + \log \lambda  + 1 + \left( {2k + 1} \right)\pi i} \right|} }}{{1 + te^{ - i\omega _k } /a}}\frac{{dt}}{{\Gamma^\ast \left( {te^{ - i\left( {\omega _k  - \pi } \right)} } \right)}}} } 
\\  =\; & \frac{{\left( { - 1} \right)^N a^N }}{{\left( {z + a} \right)^{2N + 1} }}\left( {\lambda  + 1} \right)^{2N + 1} \sum\limits_{k = 0}^\infty  {\frac{1}{{\sqrt {2\pi } i}}\int_0^{ + \infty e^{i\omega _k } } {\frac{{s^{N - \frac{1}{2}} e^{ - s\left( {\lambda  + \log \lambda  + 1 - \left( {2k + 1} \right)\pi i} \right)} }}{{1 + s/a}}\frac{{ds}}{{\Gamma^\ast \left( {se^{ - \pi i} } \right)}}} } 
\\ & + \frac{{\left( { - 1} \right)^{N + 1} a^N }}{{\left( {z + a} \right)^{2N + 1} }}\left( {\lambda  + 1} \right)^{2N + 1} \sum\limits_{k = 0}^\infty  {\frac{1}{{\sqrt {2\pi } i}}\int_0^{ + \infty e^{ - i\omega _k } } {\frac{{s^{N - \frac{1}{2}} e^{ - s\left( {\lambda  + \log \lambda  + 1 + \left( {2k + 1} \right)\pi i} \right)} }}{{1 + s/a}}\frac{{ds}}{{\Gamma^\ast \left( {se^{\pi i} } \right)}}} } .
\end{align*}
Since $\left| {\arg a} \right| <  - \omega  + \pi $ and $\frac{\pi }{2} < \omega _k  \le \omega _0  = \omega  < \pi$, we can rotate the paths of integration under the sums so that each of them has directions $\omega_0=\omega$ or $-\omega_0=-\omega$, respectively. In this way, we arrive at the representation
\begin{align*}
R_N \left( {a,\lambda } \right) = \; & \frac{{\left( { - 1} \right)^N a^N }}{{\left( {z + a} \right)^{2N + 1} }}\left( {\lambda  + 1} \right)^{2N + 1} \frac{1}{{\sqrt {2\pi } i}}\int_0^{ + \infty e^{i\omega } } {\frac{{s^{N - \frac{1}{2}} e^{ - s\left( {\lambda  + \log \lambda  + 1 - \pi i} \right)} }}{{1 + s/a}}\sum\limits_{k = 0}^\infty  {e^{2\pi iks} } \frac{{ds}}{{\Gamma^\ast \left( {se^{ - \pi i} } \right)}}} 
\\ & + \frac{{\left( { - 1} \right)^{N + 1} a^N }}{{\left( {z + a} \right)^{2N + 1} }}\left( {\lambda  + 1} \right)^{2N + 1} \frac{1}{{\sqrt {2\pi } i}}\int_0^{ + \infty e^{ - i\omega } } {\frac{{s^{N - \frac{1}{2}} e^{ - s\left( {\lambda  + \log \lambda  + 1 + \pi i} \right)} }}{{1 + s/a}}\sum\limits_{k = 0}^\infty  {e^{ - 2\pi iks} } \frac{{ds}}{{\Gamma^\ast \left( {se^{\pi i} } \right)}}} ,
\end{align*}
which holds for $\left| {\arg a} \right| <  - \omega  + \pi $, $N\geq 1$ and $0<\lambda< W\left( {e^{ - 1} } \right)$. Noting that
\[
\sum\limits_{k = 0}^\infty  {e^{ \pm 2\pi iks} } \frac{1}{{\Gamma^\ast \left( {se^{ \mp \pi i} } \right)}} = \frac{1}{{1 - e^{ \pm 2\pi iks} }}\frac{1}{{\Gamma^\ast \left( {se^{ \mp \pi i} } \right)}} = \Gamma^\ast \left( s \right),
\]
the proof of the formula \eqref{eq8} is complete.

To prove the second representation in \eqref{eq9}, we substitute \eqref{eq8} into the right hand side of
\begin{equation}\label{eq41}
b_n \left( { - \lambda } \right) = \frac{{\left( {z + a} \right)^{2n + 1} }}{{a^n }}\left( {R_n \left( {a,\lambda } \right) - R_{n + 1} \left( {a,\lambda } \right)} \right).
\end{equation}

\section{Error bounds}\label{section3}
In this section, we shall give computable bounds for the remainder term $R_N \left( {a,\lambda } \right)$ of the asymptotic series \eqref{eq10}. To estimate the remainder term, we will use the elementary result
\begin{equation}\label{eq40}
\frac{1}{{\left| {1 - re^{i\varphi } } \right|}} \le \begin{cases} \left|\csc \varphi \right| & \; \text{ if } \; 0 < \left|\varphi \text{ mod } 2\pi\right| <\frac{\pi}{2}, \\ 1 & \; \text{ if } \; \frac{\pi}{2} \leq \left|\varphi \text{ mod } 2\pi\right| \leq \pi, \end{cases}
\end{equation}
which holds for any $r>0$. Most of the error bounds we shall consider are valid in the sector $\left| {\arg a} \right| <  - \omega  + \pi$. Estimations for $R_N \left( {a,\lambda } \right)$ beyond the lines $\arg a =  \mp \omega  \pm \pi$ may be obtained using techniques similar to that we apply in Appendix \ref{appendixa} for the error term $M_N \left( z \right)$ of the asymptotic series of the scaled gamma function. We do not pursue the details here. We split the analysis into four parts depending on the size of the parameter $\lambda$.

\subsection{Case (i): $\lambda > W\left( {e^{\pi  - 1} } \right)$} If $\lambda > W\left( {e^{\pi  - 1} } \right) = 1.64428 \ldots$, then $0<\omega<\frac{\pi}{4}$. Simple estimation of \eqref{eq8}, the inequality \eqref{eq40} and the identity $\left| {\Gamma^\ast \left( {te^{i\omega } } \right)} \right| = \left| {\Gamma^\ast \left( {te^{ - i\omega } } \right)} \right|$ yields the bound
\begin{gather}\label{eq42}
\begin{split}
& \left| {R_N \left( {a,\lambda } \right)} \right| \le \left| {\frac{{a^N }}{{\left( {z + a} \right)^{2N + 1} }}} \right|\frac{{2\left( {\lambda  + 1} \right)^{2N + 1} }}{{\sqrt {2\pi } }}\int_0^{ + \infty } {t^{N - \frac{1}{2}} e^{ - t\left| {\lambda  + \log \lambda  + 1 + \pi i} \right|} \left| {\Gamma^\ast \left( {te^{i\omega } } \right)} \right|dt} 
\\ & \times\frac{1}{2}\left( \begin{cases} \left| {\csc \left( {\theta  - \omega } \right)} \right| & \text{ if } \omega  - \pi  < \theta  < \omega  - \frac{\pi }{2} \text{ or }
\omega  + \frac{\pi }{2} < \theta  <  - \omega  + \pi, \\ 1 &  \text{ if }  \omega  - \frac{\pi }{2} \le \theta  \le \omega  + \frac{\pi }{2} \end{cases}\right.
\\ & +\left.\begin{cases} \left| {\csc \left( {\theta  + \omega } \right)} \right| & \text{ if } \omega  - \pi  < \theta  < -\omega  - \frac{\pi }{2} \text{ or }
-\omega  + \frac{\pi }{2} < \theta  <  - \omega  + \pi, \\ 1 & \text{ if }  -\omega  - \frac{\pi}{2} \le \theta  \le -\omega  + \frac{\pi }{2} \end{cases}\right),
\end{split}
\end{gather}
with $\theta = \arg a$. We can simplify further this bound by estimating the quantity $\left| {\Gamma^\ast \left( {te^{i\omega } } \right)} \right|$ under the integral. Employing the notation in \eqref{eq11}, we can write
\[
\left| {\Gamma^\ast \left( {te^{i\omega } } \right)} \right| = \left| {1  + M_1 \left( {te^{i\omega } } \right)} \right| \le 1 + \left| {M_1 \left( {te^{i\omega } } \right)} \right|.
\]
The following sharp estimate was proved in \cite{Nemes}:
\[
\left| {M_1 \left( {te^{i\omega } } \right)} \right| \le \frac{1}{12t}+\frac{1}{288 t^2}
\]
for $0<\omega<\frac{\pi}{4}$, whence the quantity in the first line of \eqref{eq42} is bounded from above by
\begin{gather}\label{eq43}
\begin{split}
& \left| {\frac{{a^N }}{{\left( {z + a} \right)^{2N + 1} }}} \right|\frac{{\Gamma \left( {N + \frac{1}{2}} \right)\left( {\lambda  + 1} \right)^{2N + 1} }}{{\left( {\frac{1}{2}\pi \left| {\lambda  + \log \lambda  + 1 + \pi i} \right|} \right)^{\frac{1}{2}} \left| {\lambda  + \log \lambda  + 1 + \pi i} \right|^N }} \\ & \times \left( {1 + \frac{{\left| {\lambda  + \log \lambda  + 1 + \pi i} \right|}}{{12\left( {N - \frac{1}{2}} \right)}} + \frac{{\left| {\lambda  + \log \lambda  + 1 + \pi i} \right|^2 }}{{288\left( {N - \frac{1}{2}} \right)\left( {N - \frac{3}{2}} \right)}}} \right),
\end{split}
\end{gather}
provided that $N\geq 2$. We would like to make sure that the resulting error estimate is realistic, that is, it does not seriously overestimate the actual error. From the results in Section \ref{section4}, for large $N$, the $N$th term of the asymptotic series \eqref{eq10} satisfies
\begin{equation}\label{eq45}
\left| {\frac{{a^N b_N \left( { - \lambda } \right)}}{{\left( {z + a} \right)^{2N + 1} }}} \right| \sim \left| {\frac{{a^N }}{{\left( {z + a} \right)^{2N + 1} }}} \right|\frac{{\Gamma \left( {N + \frac{1}{2}} \right)\left( {\lambda  + 1} \right)^{2N + 1} }}{{\left( {\frac{1}{2}\pi \left| {\lambda  + \log \lambda  + 1 + \pi i} \right|} \right)^{\frac{1}{2}} \left| {\lambda  + \log \lambda  + 1 + \pi i} \right|^N }}\left| {\sin \left( {\left( {N + \frac{1}{2}} \right)\omega } \right)} \right|,
\end{equation}
as long as $\sin \left( {\left( {N + \frac{1}{2}} \right)\omega } \right) \neq 0$. Therefore if $\left|\arg a\right|$ is not very close to $-\omega+\pi$ and $N$ is large, the error bound resulting from the combination of \eqref{eq42} and \eqref{eq43} is indeed realistic.

\subsection{Case (ii): $\lambda  \ge W\left( {e^{ - 1} } \right)$} If $\lambda  \ge W\left( {e^{ - 1} } \right) = 0.27846 \ldots$, then $0<\omega\leq\frac{\pi}{2}$. We give an error bound that is valid in this $\lambda$-region and in the sector $\left|\arg a\right|\leq\frac{\pi}{4}$. If we use the relation
\[
e^{ - \left( {N + \frac{1}{2}} \right)\omega i} \Gamma^\ast \left( {te^{ - i\omega } } \right) = \overline {e^{\left( {N + \frac{1}{2}} \right)\omega i} \Gamma^\ast \left( {te^{i\omega } } \right)}
\]
in \eqref{eq10}, a straightforward calculation yields
\begin{gather}\label{eq44}
\begin{split}
& R_N \left( {a,\lambda } \right) = \sin \left( {\left( {N + \frac{1}{2}} \right)\omega } \right)\frac{{\left( { - 1} \right)^N a^N }}{{\left( {z + a} \right)^{2N + 1} }}\left( {\lambda  + 1} \right)^{2N + 1} \sqrt {\frac{2}{\pi }} \int_0^{ + \infty } {\frac{{t^{N - \frac{1}{2}} e^{ - t\left| {\lambda  + \log \lambda  + 1 + \pi i} \right|} }}{{\left( {1 + te^{i\omega } /a} \right)\left( {1 + te^{ - i\omega } /a} \right)}}\Re \Gamma^\ast \left( {te^{i\omega } } \right)dt} 
\\ & + \cos \left( {\left( {N + \frac{1}{2}} \right)\omega } \right)\frac{{\left( { - 1} \right)^N a^N }}{{\left( {z + a} \right)^{2N + 1} }}\left( {\lambda  + 1} \right)^{2N + 1} \sqrt {\frac{2}{\pi }} \int_0^{ + \infty } {\frac{{t^{N - \frac{1}{2}} e^{ - t\left| {\lambda  + \log \lambda  + 1 + \pi i} \right|} }}{{\left( {1 + te^{i\omega } /a} \right)\left( {1 + te^{ - i\omega } /a} \right)}}\Im \Gamma^\ast \left( {te^{i\omega } } \right)dt} 
\\ & + \sin \left( {\left( {N - \frac{1}{2}} \right)\omega } \right)\frac{{\left( { - 1} \right)^N a^{N + 1} }}{{\left( {z + a} \right)^{2N + 3} }}\left( {\lambda  + 1} \right)^{2N + 3} \sqrt {\frac{2}{\pi }} \int_0^{ + \infty } {\frac{{t^{N + \frac{1}{2}} e^{ - t\left| {\lambda  + \log \lambda  + 1 + \pi i} \right|} }}{{\left( {1 + te^{i\omega } /a} \right)\left( {1 + te^{ - i\omega } /a} \right)}}\Re \Gamma^\ast \left( {te^{i\omega } } \right)dt} 
\\ & + \cos \left( {\left( {N - \frac{1}{2}} \right)\omega } \right)\frac{{\left( { - 1} \right)^N a^{N + 1} }}{{\left( {z + a} \right)^{2N + 3} }}\left( {\lambda  + 1} \right)^{2N + 3} \sqrt {\frac{2}{\pi }} \int_0^{ + \infty } {\frac{{t^{N + \frac{1}{2}} e^{ - t\left| {\lambda  + \log \lambda  + 1 + \pi i} \right|} }}{{\left( {1 + te^{i\omega } /a} \right)\left( {1 + te^{ - i\omega } /a} \right)}}\Im \Gamma^\ast \left( {te^{i\omega } } \right)dt} .
\end{split}
\end{gather}
It can be shown that
\begin{equation}\label{eq46}
\frac{1}{{\left| {\left( {1 + te^{i\omega } /a} \right)\left( {1 + te^{ - i\omega } /a} \right)} \right|}} = \frac{1}{{\left| {\left( {1 + te^{i\left( {\omega  - \theta } \right)} /\left| a \right|} \right)\left( {1 + te^{ - i\left( {\omega  + \theta } \right)} /\left| a \right|} \right)} \right|}} \le 1,
\end{equation}
if $0<\omega\leq\frac{\pi}{2}$ and $\left|\arg a\right|\leq\frac{\pi}{4}$. From the results of the paper \cite{Nemes}, it follows that
\[
\left| {M_2 \left( {te^{i\omega } } \right)} \right| \le \frac{{2\sqrt 2  + 1}}{2}\frac{{\left( {1 + \zeta \left( 2 \right)} \right)\Gamma \left( 2 \right)}}{{\left( {2\pi } \right)^3 t^2 }} < \frac{1}{48t^2}
\]
for $0<\omega\leq\frac{\pi}{2}$, whence the quantities $\Re \Gamma^\ast \left( {te^{i\omega } } \right)$ and $\Im \Gamma^\ast \left( {te^{i\omega } } \right)$ can be bounded as follows
\[
\left| {\Re \Gamma^\ast \left( {te^{i\omega } } \right)} \right| = \left| {1 + \frac{{\cos \omega }}{{12t}} + \Re M_2 \left( {te^{i\omega } } \right)} \right| \le 1 + \frac{{\left| {\cos \omega } \right|}}{{12t}} + \left| {M_2 \left( {te^{i\omega } } \right)} \right| \le 1 + \frac{{\left| {\cos \omega } \right|}}{{12t}} + \frac{1}{{48t^2 }},
\]
\[
\left| {\Im \Gamma^\ast \left( {te^{i\omega } } \right)} \right| = \left| {-\frac{{\sin \omega }}{{12t}} + \Im M_2 \left( {te^{i\omega } } \right)} \right| \le \frac{{\left| {\sin \omega } \right|}}{{12t}} + \left| {M_2 \left( {te^{i\omega } } \right)} \right| \le \frac{{\left| {\sin \omega } \right|}}{{12t}} + \frac{1}{{48t^2 }}.
\]
Now a simple estimation of the representation \eqref{eq44}, the inequality \eqref{eq46}, and the bounds for $\Re \Gamma^\ast \left( {te^{i\omega } } \right)$ and $\Im \Gamma^\ast \left( {te^{i\omega } } \right)$ produce a bound for $R_N \left( {a,\lambda } \right)$ with $N\geq 2$, valid in the sector $\left|\arg a\right|\leq\frac{\pi}{4}$. This bound contains more terms than the one we have derived in Case (i), but for large $N$ and $a$ it is asymptotic to the right-hand side of \eqref{eq45}, therefore it is more realistic.

\subsection{Case (iii): $W\left( {e^{ - \pi  - 1} } \right) \le \lambda  \le W\left( {e^{\pi  - 1} } \right)$} If $0.01565\ldots=W\left( {e^{ - \pi  - 1} } \right) \le \lambda  \le W\left( {e^{\pi  - 1} } \right)= 1.64428 \ldots$, then $\frac{\pi }{4} \le \omega  \le \frac{{3\pi }}{4}$. Trivial estimation of \eqref{eq8}, the inequality \eqref{eq40} and the relation $\left| {\Gamma^\ast \left( {te^{i\omega } } \right)} \right| = \left| {\Gamma^\ast \left( {te^{ - i\omega } } \right)} \right|$ gives the bound
\begin{align*}
& \left| {R_N \left( {a,\lambda } \right)} \right| \le \left| {\frac{{a^N }}{{\left( {z + a} \right)^{2N + 1} }}} \right|\frac{{2\left( {\lambda  + 1} \right)^{2N + 1} }}{{\sqrt {2\pi } }}\int_0^{ + \infty } {t^{N - \frac{1}{2}} e^{ - t\left| {\lambda  + \log \lambda  + 1 + \pi i} \right|} \left| {\Gamma^\ast \left( {te^{i\omega } } \right)} \right|dt} 
\\ & \times\frac{1}{2}\left( \begin{cases} \left| {\csc \left( {\theta  - \omega } \right)} \right| & \text{ if } \omega  - \pi  < \theta  < \omega  - \frac{\pi }{2}, \\ 1 &  \text{ if }  \omega  - \frac{\pi }{2} \le \theta  \le -\omega  + \pi \end{cases} +\begin{cases} \left| {\csc \left( {\theta  + \omega } \right)} \right| & \text{ if } -\omega  + \frac{\pi }{2} < \theta  <  - \omega  + \pi, \\ 1 & \text{ if }  \omega  - \pi \le \theta  \le -\omega  + \frac{\pi }{2} \end{cases}\right).
\end{align*}
We can again simplify our estimate by introducing an upper bound for $\left| {\Gamma^\ast \left( {te^{i\omega } } \right)} \right|$, by applying the results of Appendix \ref{appendixa} to $M_2\left(te^{i\omega}\right)$. The resulting estimate is again comparable with the quantity on the right-hand side of \eqref{eq45} if $\left|\arg a\right|$ is not very close to $-\omega+\pi$ and $N$ is large.

\subsection{Case (iv): $0< \lambda  < W\left( {e^{-\pi  - 1} } \right)$} If $0< \lambda  \le W\left( {e^{-\pi  - 1} } \right)=0.01565 \ldots$, then $\frac{3\pi }{4} <\omega <\pi$. The method of Case (i) and (iii) still applies and yields
\begin{align*}
\left| {R_N \left( {a,\lambda } \right)} \right| \le \; &\left| {\frac{{a^N }}{{\left( {z + a} \right)^{2N + 1} }}} \right|\frac{{2\left( {\lambda  + 1} \right)^{2N + 1} }}{{\sqrt {2\pi } }}\int_0^{ + \infty } {t^{N - \frac{1}{2}} e^{ - t\left| {\lambda  + \log \lambda  + 1 + \pi i} \right|} \left| {\Gamma^\ast \left( {te^{i\omega } } \right)} \right|dt} 
\\ & \times \frac{{\left| {\csc \left( {\theta  - \omega } \right)} \right| + \left| {\csc \left( {\theta  + \omega } \right)} \right|}}{2},
\end{align*}
as long as $\left| {\arg a} \right| <  - \omega  + \pi$. This bound, however, is not very effective when $\lambda$ is close to $0$ or, equivalently, $\omega$ is close to $\pi$. When $\lambda$ approaches $0$, $te^{i\omega }$ ($t>0$) approaches the negative real axis where the poles of the gamma function lie; moreover the sector of validity becomes smaller and smaller due to the poles of the cosecant at $\pm \pi$.

\section{Asymptotics for the late coefficients}\label{section4} In this section, we investigate the asymptotic nature of the coefficients $b_n \left( -\lambda  \right)$ as $n\to +\infty$. For our purposes, the most appropriate representation of these coefficients is the second formula in \eqref{eq9}. Upon replacing $\Gamma^\ast \left( {te^{ \pm i\omega } } \right)$ by their representation \eqref{eq11} in this integral, we obtain
\begin{gather}\label{eq12}
\begin{split}
b_n \left( { - \lambda } \right) = \; & \frac{{\left( { - 1} \right)^n \Gamma \left( {n + \frac{1}{2}} \right)\left( {\lambda  + 1} \right)^{2n + 1} }}{{\left( {\frac{1}{2}\pi \left| {\lambda  + \log \lambda  + 1 + \pi i} \right|} \right)^{\frac{1}{2}} \left| {\lambda  + \log \lambda  + 1 + \pi i} \right|^n }} \\ & \times \left( {\sum\limits_{k = 0}^{K - 1} {\left( { - 1} \right)^k \left| {\lambda  + \log \lambda  + 1 + \pi i} \right|^k \gamma _k \frac{{\Gamma \left( {n - k + \frac{1}{2}} \right)}}{{\Gamma \left( {n + \frac{1}{2}} \right)}}\sin \left( {\left( {n - k + \frac{1}{2}} \right)\omega } \right)}  + A_K \left( {n,\lambda } \right)} \right),
\end{split}
\end{gather}
for any fixed $1\leq K \leq n-1$, provided that $n\geq 2$. The remainder term $A_K \left( {n,\lambda } \right)$ is given by the integral formula
\begin{align*}
A_K \left( {n,\lambda } \right) = \frac{{\left| {\lambda  + \log \lambda  + 1 + \pi i} \right|^{n + \frac{1}{2}} }}{{2i\Gamma \left( {n + \frac{1}{2}} \right)}} & \int_0^{ + \infty } t^{n - \frac{1}{2}} e^{ - t\left| {\lambda  + \log \lambda  + 1 + \pi i} \right|} \\ & \times \left( {e^{\left( {n + \frac{1}{2}} \right)\omega i} M_K \left( {te^{i\omega } } \right) - e^{ - \left( {n + \frac{1}{2}} \right)\omega i} M_K \left( {te^{ - i\omega } } \right)} \right)dt .
\end{align*}
Since $\left| {M_K \left( {te^{i\omega } } \right)} \right| = \left| {\overline {M_K \left( {te^{i\omega } } \right)} } \right| = \left| {M_K \left( {\overline {te^{i\omega } } } \right)} \right| = \left| {M_K \left( {te^{ - i\omega } } \right)} \right|$, trivial estimation yields the bound
\[
\left| {A_K \left( {n,\lambda } \right)} \right| \le \frac{{\left| {\lambda  + \log \lambda  + 1 + \pi i} \right|^{n + \frac{1}{2}} }}{{\Gamma \left( {n + \frac{1}{2}} \right)}}\int_0^{ + \infty } {t^{n - \frac{1}{2}} e^{ - t\left| {\lambda  + \log \lambda  + 1 + \pi i} \right|} \left| {M_K \left( {te^{i\omega } } \right)} \right|dt} .
\]
First assume that $\lambda  \ge W\left( {e^{\pi  - 1} } \right) = 1.64428 \ldots$, i.e., $0 < \omega  \le \frac{\pi }{4}$. In this case, it was proved in \cite{Nemes} that
\[
\left| {M_K \left( {te^{i\omega } } \right)} \right| \le \frac{{\left| {\gamma _K } \right|}}{{t^K }} + \frac{{\left| {\gamma _{K + 1} } \right|}}{{t^{K + 1} }}
\]
for any $t>0$, which leads to the simple estimate
\begin{gather}\label{eq17}
\begin{split}
\left| {A_K \left( {n,\lambda } \right)} \right| \le \; & \left| {\lambda  + \log \lambda  + 1 + \pi i} \right|^K \left| {\gamma _K } \right|\frac{{\Gamma \left( {n - K + \frac{1}{2}} \right)}}{{\Gamma \left( {n + \frac{1}{2}} \right)}} \\ &+ \left| {\lambda  + \log \lambda  + 1 + \pi i} \right|^{K + 1} \left| {\gamma _{K + 1} } \right|\frac{{\Gamma \left( {n - K - \frac{1}{2}} \right)}}{{\Gamma \left( {n + \frac{1}{2}} \right)}}.
\end{split}
\end{gather}
For the general case $\lambda>0$, i.e., when $0<\omega <\pi$, we prove in Appendix \ref{appendixa} that
\[
\left| {M_K \left( {te^{i\omega } } \right)} \right| \le \left( {\frac{{\sec \left( {\omega  - \varphi^\ast } \right)}}{{\cos ^K \varphi^\ast }} + 1} \right)\frac{{\zeta \left( K \right)\Gamma \left( K \right)}}{{\left( {2\pi } \right)^{K + 1} t^K }}
\]
for any $t>0$ and $K\geq 2$, where $\varphi^\ast$ is the unique solution of the equation
\[
\left( {K + 1} \right)\sin \left( {\omega  - 2\varphi^\ast } \right) = \left( {K - 1} \right)\sin \omega ,
\]
that satisfies $ - \frac{\pi}{2} + \omega  < \varphi ^\ast < \frac{\pi }{2} $ if $ \frac{\pi}{2} \le \omega  < \pi$, and $0 < \varphi^\ast  < \omega$ if $0  < \omega  < \frac{\pi}{2}$. With this estimate, we obtain the error bound
\begin{equation}\label{eq13}
\left| {A_K \left( {n,\lambda } \right)} \right| \le \frac{1}{2}\left( {\frac{{\sec \left( {\omega  - \varphi^\ast } \right)}}{{\cos ^K \varphi^\ast }} + 1} \right)\left| {\lambda  + \log \lambda  + 1 + \pi i} \right|^K \frac{1}{\pi}\frac{{\zeta \left( K \right)\Gamma \left( K \right)}}{{\left( {2\pi } \right)^{K} }}\frac{{\Gamma \left( {n - K + \frac{1}{2}} \right)}}{{\Gamma \left( {n + \frac{1}{2}} \right)}}.
\end{equation}
We remark that if $K$ is odd and large, $\frac{1}{\pi} \frac{{\zeta \left( K \right)\Gamma \left( K \right)}}{{\left( {2\pi } \right)^{K} }} \sim \left| {\gamma _K } \right|$ (see, e.g., Boyd \cite{Boyd}), whence the form of this bound is closely related to the first omitted term of the series \eqref{eq12}. If $\lambda \to 0+$, or equivalently $\omega \to \pi-0$, then $\varphi^\ast \to \frac{\pi}{2}$ and the bound \eqref{eq13} becomes singular. This singular behaviour of the error bound is related to the poles of the gamma function along the negative real axis. Whence, the expansion \eqref{eq12} is effective only if $\lambda$ is not very close to $0$. An alternative integral representation for the coefficients $b_n \left( -\lambda  \right)$, when $\lambda$ is close to $0$, can be obtained by combining Theorem \ref{thm4} with \eqref{eq41}. For the resulting expression one may use the truncated version of the asymptotic series of the reciprocal gamma function together with the known sharp error bounds \cite{Nemes}. The details are not discussed here.

If $\lambda$ is fixed and is not very close to $0$ and $n$ is large, the least value of the bound \eqref{eq13} occurs when
\[
K \approx \left( {n + \frac{1}{2}} \right)\frac{{2\pi }}{{\left| {\lambda  + \log \lambda  + 1 + \pi i} \right| + 2\pi }}.
\]
With this choice of $K$, the error bound is
\[
\mathcal{O}\left( {n^{ - \frac{1}{2}} \left( {\frac{{\left| {\lambda  + \log \lambda  + 1 + \pi i} \right|}}{{\left| {\lambda  + \log \lambda  + 1 + \pi i} \right| + 2\pi }}} \right)^n } \right).
\]
This is the best accuracy we can achieve using the expansion \eqref{eq12}. Whence, the larger $\lambda$ is the larger $n$ has to be to get a reasonable approximation from \eqref{eq12}.

\begin{table*}[!ht]
\begin{center}
\begin{tabular}
[c]{ l r @{\,}c@{\,} l}\hline
 & \\ [-1ex]
 values of $\lambda$ and $K$ & $\lambda=\frac{1}{100}$, $K=57$ & & \\ [1ex]
 exact numerical value of $b_{100}\left(-\lambda\right)$ & $-0.320681358577665454823220737555930836965007363$ & $\times$ & $10^{90}$ \\ [1ex]
 approximation \eqref{eq12} to $b_{100}\left(-\lambda\right)$ & $-0.320681358577665454823220737555930624925282126$ & $\times$ & $10^{90}$  \\ [1ex]
 error & $-0.212039725237$ & $\times$ & $10^{57}$\\ [1ex] 
 error bound using \eqref{eq13} & $0.23677013448560065229$ & $\times$ & $10^{65}$\\ [1ex] \hline
 & \\ [-1ex]
 values of $\lambda$ and $K$ & $\lambda=2$, $K=57$ & & \\ [1ex]
 exact numerical value of $b_{100}\left(-\lambda\right)$ & $0.732252465623483776580694573188048575042344276$ & $\times$ & $10^{184}$ \\ [1ex]
 approximation \eqref{eq12} to $b_{100}\left(-\lambda\right)$ & $0.732252465623483776580694573188048575045373691$ & $\times$ & $10^{184}$ \\ [1ex]
 error & $-0.3029415$ & $\times$ & $10^{146}$\\ [1ex]
 error bound using \eqref{eq17} & $0.63782498$ & $\times$ & $10^{147}$\\ [1ex] \hline
 & \\ [-1ex]
 values of $\lambda$ and $K$ &  $\lambda=5$, $K=43$ & & \\ [1ex]
 exact numerical value of $b_{100}\left(-\lambda\right)$ & $0.186478888380183206402841100236655575383457561$ & $\times$ & $10^{222}$ \\ [1ex]
 approximation \eqref{eq12} to $b_{100}\left(-\lambda\right)$ & $0.186478888380183206402841097953515994081833820$ & $\times$ & $10^{222}$ \\ [1ex]
 error & $0.2283139581301623741$ & $\times$ & $10^{196}$\\ [1ex] 
 error bound using \eqref{eq17} & $0.5373934537697861855$ & $\times$ & $10^{196}$\\ [-1ex]
 & \\\hline
\end{tabular}
\end{center}
\caption{Approximations for $b_{100}\left(-\lambda\right)$ with various $\lambda$, using \eqref{eq12}.}
\label{table1}
\end{table*}

By extending the sum in \eqref{eq12} to infinity, we arrive at the formal series
\begin{align*}
\left( { - 1} \right)^{n + 1} & \frac{{b_n \left( { - \lambda } \right)}}{{\lambda \left( {\lambda  + 1} \right)^{2n} }} \approx  - \frac{{\Gamma \left( {n + \frac{1}{2}} \right)\left( {1 + \frac{1}{\lambda }} \right)}}{{\left( {\frac{1}{2}\pi \left| {\lambda  + \log \lambda  + 1 + \pi i} \right|} \right)^{\frac{1}{2}} \left| {\lambda  + \log \lambda  + 1 + \pi i} \right|^n }}\left( {\sin \left( {\left( {n + \frac{1}{2}} \right)\omega } \right) }\right.  \\ & \left.{+ \frac{{\left| {\lambda  + \log \lambda  + 1 + \pi i} \right|}}{{12\left( {n - \frac{1}{2}} \right)}}\sin \left( {\left( {n - \frac{1}{2}} \right)\omega } \right) + \frac{{\left| {\lambda  + \log \lambda  + 1 + \pi i} \right|^2 }}{{288\left( {n - \frac{1}{2}} \right)\left( {n - \frac{3}{2}} \right)}}\sin \left( {\left( {n - \frac{3}{2}} \right)\omega } \right) }\right.  \\ & \left.{- \frac{{139\left| {\lambda  + \log \lambda  + 1 + \pi i} \right|^3 }}{{51840\left( {n - \frac{1}{2}} \right)\left( {n - \frac{3}{2}} \right)\left( {n - \frac{5}{2}} \right)}}\sin \left( {\left( {n - \frac{5}{2}} \right)\omega } \right) +  \cdots } \right).
\end{align*}
This is exactly Dingle's expansion for the late coefficients in the asymptotic series of $\Gamma\left(-a,z\right)$ \cite[p. 162]{Dingle}. Note that Dingle denotes $\omega$ by $\theta$ and defines it as
\[
\arctan \left( {\frac{\pi }{{\lambda  + \log \lambda  + 1}}} \right).
\]
This expression for $\omega$ is, however, correct only when $\lambda > W\left( {e^{ - 1} } \right) =0.27846\ldots$, or equivalently when $0 < \omega < \frac{\pi}{2}$.  The mathematically rigorous form of Dingle's series is therefore the formula \eqref{eq12}.

Numerical examples illustrating the efficacy of the expansion \eqref{eq12}, truncated optimally, are given in Table \ref{table1}.

\section{Exponentially improved asymptotic expansion}\label{section5}

The aim of this section is to provide a rigorous treatment of Dingle's formal re-expansion of the remainder term $R_N \left( {a,\lambda } \right)$ of the asymptotic series \eqref{eq10} \cite[p. 463]{Dingle}. The main result is stated in Theorem \ref{thm2} below. In this theorem, we truncate the asymptotic series \eqref{eq10} at about its least term and re-expand the remainder into a new asymptotic expansion. The resulting exponentially improved asymptotic series is valid in a larger region than the original expansion. The terms in this new series involve the terminant function $\widehat T_p\left(w\right)$ (see below), which allows the smooth transition through the Stokes lines. Throughout this section, we use subscripts in the $\mathcal{O}$ notations to indicate the dependence of the implied constant on certain parameters. In this theorem, $R_N \left( {a,\lambda } \right)$ is extended to a wider region using analytic continuation.

\begin{theorem}\label{thm2} Let $K$ be an arbitrary fixed non-negative integer, and let $\lambda>0$ be a fixed real number. Suppose that $\left| {\arg a} \right| \le \max \left( {\omega  + \pi ,2\pi  - \delta } \right)<2\pi$ with an arbitrary fixed small positive $\delta$, $\left|a\right|$ is large and $N = \left| a \left( {\lambda  + \log \lambda  + 1+\pi i} \right)\right|+ \rho$ with $\rho$ being bounded. Then
\begin{gather}\label{eq29}
\begin{split}
R_N \left( {a,\lambda } \right) = & - ie^{a\left( {\lambda  + \log \lambda  + 1 - \pi i} \right)} \sqrt {\frac{{2\pi }}{a}} \sum\limits_{k = 0}^{K - 1} {\frac{{\gamma _k }}{{a^k }}\widehat T_{N - k + \frac{1}{2}} \left( {a\left( {\lambda  + \log \lambda  + 1 - \pi i} \right)} \right)} 
\\ & + ie^{a\left( {\lambda  + \log \lambda  + 1 + \pi i} \right)} \sqrt {\frac{{2\pi }}{a}} \sum\limits_{k = 0}^{K - 1} {\frac{{\gamma _k }}{{a^k }}\widehat T_{N - k + \frac{1}{2}} \left( {a\left( {\lambda  + \log \lambda  + 1 + \pi i} \right)} \right)} 
\\ & + R_{N,K} \left( {a,\lambda } \right),
\end{split}
\end{gather}
where
\[
R_{N,K} \left( {a,\lambda } \right) = \mathcal{O}_{\lambda ,K,\rho } \left( {\frac{{e^{ - \left| {a\left( {\lambda  + \log \lambda  + 1 + \pi i} \right)} \right|} }}{{\left| a \right|^{K + \frac{1}{2}} }}} \right)
\]
for $\left| {\arg a} \right| \leq - \omega  + \pi$;
\[
R_{N,K} \left( {a,\lambda } \right) = \mathcal{O}_{\lambda ,K,\rho ,\delta } \left( {\frac{{e^{\left| {a\left( {\lambda  + \log \lambda  + 1 + \pi i} \right)} \right|\cos \left( {\arg a \pm \omega } \right)} }}{{\left| a \right|^{K + \frac{1}{2}} }}} \right)
\]
for $ - \omega  + \pi  \le  \pm \arg a \le \max \left( {\omega  + \pi ,2\pi  - \delta } \right)$.
\end{theorem}

It has to be noted that the expansion \eqref{eq29} is only of theoretical interest, since the terminant functions on the right-hand side are in general more complicated functions than the incomplete gamma function $\Gamma\left(-a,z\right)$, we want to approximate. If the sum in \eqref{eq29} is extended formally to infinity and the error term $R_{N,K} \left( {a,\lambda } \right)$ is neglected, the result is equivalent to Dingle's formal expansion. We remark that Dingle himself did not consider the region of validity of his expansion.

While proving Theorem \ref{thm2}, we also obtain the following explicit bound for the remainder in \eqref{eq29}. Note that in this theorem $N$ may not depend on $a$ and $\lambda$.

\begin{theorem}\label{thm3} For any integers $2 \leq K \leq N$, define the remainder $R_{N,K} \left( {a,\lambda } \right)$ by \eqref{eq29}. Then we have
\begin{gather}\label{eq37}
\begin{split}
& \left| R_{N,K}\left(a,\lambda\right) \right|
\\ & \le \left( {\frac{{\sec \left( {\omega  - \varphi^\ast } \right)}}{{\cos ^K \varphi^\ast }} + 1} \right)\left| {e^{a\left( {\lambda  + \log \lambda  + 1 - \pi i} \right)} \sqrt {\frac{{2\pi }}{a}} \widehat T_{N - K + \frac{1}{2}} \left( {a\left( {\lambda  + \log \lambda  + 1 - \pi i} \right)} \right)} \right|\frac{{\zeta \left( K \right)\Gamma \left( K \right)}}{{\left( {2\pi } \right)^{K + 1} \left| a \right|^K }}
\\ & + \left( {\frac{{\sec \left( {\omega  - \varphi^\ast } \right)}}{{\cos ^K \varphi^\ast }} + 1} \right)\left| {e^{a\left( {\lambda  + \log \lambda  + 1 + \pi i} \right)} \sqrt {\frac{{2\pi }}{a}} \widehat T_{N - K + \frac{1}{2}} \left( {a\left( {\lambda  + \log \lambda  + 1 + \pi i} \right)} \right)} \right|\frac{{\zeta \left( K \right)\Gamma \left( K \right)}}{{\left( {2\pi } \right)^{K + 1} \left| a \right|^K }}
\\ & + \left( {\frac{{\sec ^2 \left( {\omega  - \varphi^\ast } \right)}}{{\cos ^K \varphi^\ast }} + 1} \right)\frac{{2\zeta \left( K \right)\Gamma \left( K \right)\Gamma \left( {N - K + \frac{1}{2}} \right)}}{{\left( {2\pi } \right)^{K + \frac{3}{2}} \left| {\lambda  + \log \lambda  + 1 + \pi i} \right|^{N - K + \frac{1}{2}} \left| a \right|^{N + 1} }},
\end{split}
\end{gather}
provided that $\left| {\arg a} \right| \leq - \omega  + \pi$. Here $0< \varphi^\ast <\frac{\pi}{2}$ is the unique solution of the equation
\[
\left( {K + 1} \right)\sin \left( {\omega  - 2\varphi^\ast } \right) = \left( {K - 1} \right)\sin \omega ,
\]
that satisfies $ - \frac{\pi}{2} + \omega  < \varphi ^\ast < \frac{\pi }{2} $ if $ \frac{\pi}{2} \le \omega < \pi$, and $0 < \varphi^\ast  < \omega$ if $0  < \omega  < \frac{\pi}{2}$.
\end{theorem}

The (scaled) terminant function can be defined in terms of the incomplete gamma function as
\[
\widehat T_p \left( w \right) = \frac{{e^{\pi ip} \Gamma \left( p \right)}}{{2\pi i}}\Gamma \left( {1 - p,w} \right) = \frac{e^{\pi ip} w^{1 - p} e^{ - w} }{2\pi i}\int_0^{ + \infty } {\frac{{t^{p - 1} e^{ - t} }}{w + t}dt} \; \text{ for } \; p>0 \; \text{ and } \; \left| \arg w \right| < \pi ,
\]
and by analytic continuation elsewhere. Olver \cite[equations (4.5) and (4.6)]{Olver4} showed that when $p \sim \left|w\right|$ and $w \to \infty$, we have
\begin{equation}\label{eq22}
\widehat T_p \left( w \right) = \begin{cases} \mathcal{O}\left( {e^{ - w - \left| w \right|} } \right) & \; \text{ if } \; \left| {\arg w} \right| \le \pi, \\ \mathcal{O}\left(1\right) & \; \text{ if } \; - 3\pi  < \arg w \le  - \pi. \end{cases}
\end{equation}
Concerning the smooth transition of the Stokes discontinuities, one may use the more precise asymptotic formulas
\begin{equation}\label{eq23}
\widehat T_p \left( w \right) = \frac{1}{2} + \frac{1}{2}\mathop{\text{erf}} \left( {c\left( \varphi  \right)\sqrt {\frac{1}{2}\left| w \right|} } \right) + \mathcal{O}\left( {\frac{{e^{ - \frac{1}{2}\left| w \right|c^2 \left( \varphi  \right)} }}{{\left| w \right|^{\frac{1}{2}} }}} \right)
\end{equation}
for $-\pi +\delta \leq \arg w \leq 3 \pi -\delta$, $0 < \delta  \le 2\pi$; and
\begin{equation}\label{eq24}
e^{ - 2\pi ip} \widehat T_p \left( w \right) =  - \frac{1}{2} + \frac{1}{2}\mathop{\text{erf}} \left( { - \overline {c\left( { - \varphi } \right)} \sqrt {\frac{1}{2}\left| w \right|} } \right) + \mathcal{O}\left( {\frac{{e^{ - \frac{1}{2}\left| w \right|\overline {c^2 \left( { - \varphi } \right)} } }}{{\left| w \right|^{\frac{1}{2}} }}} \right)
\end{equation}
for $- 3\pi  + \delta  \le \arg w \le \pi  - \delta$, $0 < \delta \le 2\pi$. Here $\varphi = \arg w$ and $\mathop{\text{erf}}$ denotes the error function. The quantity $c\left( \varphi  \right)$ is defined implicitly by the equation
\[
\frac{1}{2}c^2 \left( \varphi  \right) = 1 + i\left( {\varphi  - \pi } \right) - e^{i\left( {\varphi  - \pi } \right)},
\]
and corresponds to the branch of $c\left( \varphi  \right)$ which has the following expansion in the neighbourhood of $\varphi = \pi$:
\begin{equation}\label{eq25}
c\left( \varphi  \right) = \left( {\varphi  - \pi } \right) + \frac{i}{6}\left( {\varphi  - \pi } \right)^2  - \frac{1}{{36}}\left( {\varphi  - \pi } \right)^3  - \frac{i}{{270}}\left( {\varphi  - \pi } \right)^4  +  \cdots .
\end{equation}
For complete asymptotic expansions, see Olver \cite{Olver5}. We remark that Olver uses the different notation $F_p \left( w \right) = ie^{ - \pi ip} \widehat T_p \left( w \right)$ for the terminant function and the other branch of the function $c\left( \varphi  \right)$. For further properties of the terminant function, see, for example, Paris and Kaminski \cite[Chapter 6]{Paris2}.

Now we begin the proofs of Theorems \ref{thm2} and \ref{thm3}. First, we suppose that $\left| {\arg a} \right| <  - \omega  + \pi$. Let $K \geq 2$ be a fixed integer. Substituting the expression in \eqref{eq11} into \eqref{eq8} and using the definition of the terminant function we find that
\begin{align*}
R_N \left( {a,\lambda } \right) = & - ie^{a\left( {\lambda  + \log \lambda  + 1 - \pi i} \right)} \sqrt {\frac{{2\pi }}{a}} \sum\limits_{k = 0}^{K - 1} {\frac{{\gamma _k }}{{a^k }}\widehat T_{N - k + \frac{1}{2}} \left( {a\left( {\lambda  + \log \lambda  + 1 - \pi i} \right)} \right)} 
\\ & + ie^{a\left( {\lambda  + \log \lambda  + 1 + \pi i} \right)} \sqrt {\frac{{2\pi }}{a}} \sum\limits_{k = 0}^{K - 1} {\frac{{\gamma _k }}{{a^k }}\widehat T_{N - k + \frac{1}{2}} \left( {a\left( {\lambda  + \log \lambda  + 1 + \pi i} \right)} \right)} 
\\ & + R_{N,K} \left( {a,\lambda } \right),
\end{align*}
with
\begin{gather}\label{eq18}
\begin{split}
R_{N,K} \left( {a,\lambda } \right) = \; & \frac{{\left( { - 1} \right)^N }}{{a^{N + 1} }}\frac{{e^{\left( {N + \frac{1}{2}} \right)\omega i} }}{{\sqrt {2\pi } i}}\int_0^{ + \infty } {\frac{{t^{N - \frac{1}{2}} e^{ - t\left| {\lambda  + \log \lambda  + 1 + \pi i} \right|} }}{{1 + te^{i\omega } /a}}M_K \left( {te^{i\omega } } \right)dt} 
\\ & + \frac{{\left( { - 1} \right)^{N + 1} }}{{a^{N + 1} }}\frac{{e^{ - \left( {N + \frac{1}{2}} \right)\omega i} }}{{\sqrt {2\pi } i}}\int_0^{ + \infty } {\frac{{t^{N - \frac{1}{2}} e^{ - t\left| {\lambda  + \log \lambda  + 1 + \pi i} \right|} }}{{1 + te^{ - i\omega } /a}}M_K \left( {te^{ - i\omega } } \right)dt} 
\\ = \; &\left( { - 1} \right)^N \frac{{e^{\left( {N + \frac{1}{2}} \right)\left( {\omega  - \theta } \right)i} }}{{\sqrt {2\pi a} i}}\int_0^{ + \infty } {\frac{{\tau ^{N - \frac{1}{2}} e^{ - r\tau \left| {\lambda  + \log \lambda  + 1 + \pi i} \right|} }}{{1 + \tau e^{i\left( {\omega  - \theta } \right)} }}M_K \left( {r\tau e^{i\omega } } \right)d\tau } 
\\ & + \left( { - 1} \right)^{N + 1} \frac{{e^{\left( {N + \frac{1}{2}} \right)\left( { - \omega  - \theta } \right)i} }}{{\sqrt {2\pi a} i}}\int_0^{ + \infty } {\frac{{\tau ^{N - \frac{1}{2}} e^{ - r\tau \left| {\lambda  + \log \lambda  + 1 + \pi i} \right|} }}{{1 + \tau e^{i\left( { - \omega  - \theta } \right)} }}M_K \left( {r\tau e^{ - i\omega } } \right)d\tau } ,
\end{split}
\end{gather}
under the assumption that $K\leq N$. Here we have taken $a=re^{i\theta}$. We consider the estimation of the first integral after the second equality in \eqref{eq18}. By \eqref{eq19}, we have
\begin{align*}
& M_K \left( {r\tau e^{i\omega } } \right) = \frac{1}{{2\pi i}}\frac{{i^K }}{{\left( {r\tau e^{i\omega } } \right)^K }}\left( {\frac{{e^{i\varphi } }}{{\cos \varphi }}} \right)^K \int_0^{ + \infty } {\frac{{s^{K - 1} e^{ - 2\pi \frac{{se^{i\varphi } }}{{\cos \varphi }}} \Gamma^\ast \left( {\frac{{ise^{i\varphi } }}{{\cos \varphi }}} \right)}}{{1 - ise^{i\varphi } /\left( {r\tau e^{i\omega } \cos \varphi } \right)}}ds} 
\\ & - \frac{1}{{2\pi i}}\frac{{\left( { - i} \right)^K }}{{\left( {r\tau e^{i\omega } } \right)^K }}\int_0^{ + \infty } {\frac{{t^{K - 1} e^{ - 2\pi t} \Gamma^\ast \left( { - it} \right)}}{{1 + it/\left( {r\tau e^{i\omega } } \right)}}dt} 
\\ & = \frac{1}{{2\pi i}}\frac{{i^K }}{{\left( {r\tau e^{i\omega } } \right)^K }}\left( {\frac{{e^{i\varphi } }}{{\cos \varphi }}} \right)^K \int_0^{ + \infty } {\frac{{s^{K - 1} e^{ - 2\pi \frac{{se^{i\varphi } }}{{\cos \varphi }}} \Gamma^\ast \left( {\frac{{ise^{i\varphi } }}{{\cos \varphi }}} \right)}}{{1 - ise^{i\varphi } /\left( {re^{i\omega } \cos \varphi } \right)}}ds} 
\\ & + \frac{1}{{2\pi i}}\frac{{i^K }}{{\left( {r\tau e^{i\omega } } \right)^K }}\left( {\frac{{e^{i\varphi } }}{{\cos \varphi }}} \right)^K \left( {\tau  - 1} \right)\int_0^{ + \infty } {\frac{{s^{K - 1} e^{ - 2\pi \frac{{se^{i\varphi } }}{{\cos \varphi }}} \Gamma^\ast  \left( {\frac{{ise^{i\varphi } }}{{\cos \varphi }}} \right)}}{{\left( {1 + ir\tau e^{i\omega } \cos \varphi /\left( {se^{i\varphi } } \right)} \right)\left( {1 - ise^{i\varphi } /\left( {re^{i\omega } \cos \varphi } \right)} \right)}}ds} 
\\ & - \frac{1}{{2\pi i}}\frac{{\left( { - i} \right)^K }}{{\left( {r\tau e^{i\omega } } \right)^K }}\left( {\int_0^{ + \infty } {\frac{{t^{K - 1} e^{ - 2\pi t} \Gamma^\ast \left( { - it} \right)}}{{1 + it/\left( {re^{i\omega } } \right)}}dt}  + \left( {\tau  - 1} \right)\int_0^{ + \infty } {\frac{{t^{K - 1} e^{ - 2\pi t} \Gamma^\ast \left( { - it} \right)}}{{\left( {1 + r\tau e^{i\omega } /it} \right)\left( {1 + it/\left( {re^{i\omega } } \right)} \right)}}dt} } \right),
\end{align*}
with a suitable $0<\varphi<\frac{\pi}{2}$ satisfying $\varphi <\omega<\frac{\pi}{2}+\varphi$. Substitution into the first integral in \eqref{eq18} and trivial estimation yield
\begin{align*}
& \left| {\left( { - 1} \right)^N \frac{{e^{\left( {N + \frac{1}{2}} \right)\left( {\omega  - \theta } \right)i} }}{{\sqrt {2\pi a} i}}\int_0^{ + \infty } {\frac{{\tau ^{N - \frac{1}{2}} e^{ - r\tau \left| {\lambda  + \log \lambda  + 1 + \pi i} \right|} }}{{1 + \tau e^{i\left( {\omega  - \theta } \right)} }}M_K \left( {r\tau e^{i\omega } } \right)d\tau } } \right| \\ & \le \frac{1}{{\cos ^K \varphi \sqrt {2\pi r} }}\left| {\int_0^{ + \infty } {\frac{{\tau ^{N - K - \frac{1}{2}} e^{ - r\tau \left| {\lambda  + \log \lambda  + 1 + \pi i} \right|} }}{{1 + \tau e^{i\left( {\omega  - \theta } \right)} }}d\tau } } \right|\frac{1}{{2\pi r^K }}\int_0^{ + \infty } {\frac{{s^{K - 1} e^{ - 2\pi s} \left| {\Gamma^\ast \left( {\frac{{ise^{i\varphi } }}{{\cos \varphi }}} \right)} \right|}}{{\left| {1 - ise^{i\varphi } /\left( {re^{i\omega } \cos \varphi } \right)} \right|}}ds} 
\\ & + \frac{1}{{\cos ^K \varphi \sqrt {2\pi r} }}\int_0^{ + \infty } \tau ^{N - K - \frac{1}{2}} e^{ - r\tau \left| {\lambda  + \log \lambda  + 1 + \pi i} \right|} \left| {\frac{{\tau  - 1}}{{\tau  + e^{ - i\left( {\omega  - \theta } \right)} }}} \right|\\ & \times \frac{1}{{2\pi r^K }}\int_0^{ + \infty } {\frac{{s^{K - 1} e^{ - 2\pi s} \left| {\Gamma^\ast \left( {\frac{{ise^{i\varphi } }}{{\cos \varphi }}} \right)} \right|}}{{\left|\left( {1 + ir\tau e^{i\omega } \cos \varphi /\left( {se^{i\varphi } } \right)} \right)\left( {1 - ise^{i\varphi } /\left( {re^{i\omega } \cos \varphi } \right)} \right)\right|}}ds} d\tau  
\\ & + \frac{1}{{\sqrt {2\pi r} }}\left| {\int_0^{ + \infty } {\frac{{\tau ^{N - K - \frac{1}{2}} e^{ - r\tau \left| {\lambda  + \log \lambda  + 1 + \pi i} \right|} }}{{1 + \tau e^{i\left( {\omega  - \theta } \right)} }}d\tau } } \right|\frac{1}{{2\pi r^K }}\int_0^{ + \infty } {\frac{{t^{K - 1} e^{ - 2\pi t} \left| {\Gamma^\ast \left( { - it} \right)} \right|}}{{\left| {1 + it/\left( {re^{i\omega } } \right)} \right|}}dt}  \\ & + \frac{1}{{\sqrt {2\pi r} }}\int_0^{ + \infty } {\tau ^{N - K - \frac{1}{2}} e^{ - r\tau \left| {\lambda  + \log \lambda  + 1 + \pi i} \right|} \left| {\frac{{\tau  - 1}}{{\tau  + e^{ - i\left( {\omega  - \theta } \right)} }}} \right|\frac{1}{{2\pi r^K }}\int_0^{ + \infty } {\frac{{t^{K - 1} e^{ - 2\pi t} \left| {\Gamma^\ast \left( { - it} \right)} \right|}}{{\left|\left( {1 + r\tau e^{i\omega } /it} \right)\left( {1 + it/\left( {re^{i\omega } } \right)} \right)\right|}}dt} d\tau } .
\end{align*}
Noting that
\[
\left| {\frac{{\tau  - 1}}{{\tau  + e^{ - i\left( {\omega  - \theta } \right)} }}} \right| \le 1,\; \frac{1}{{\left| {1 + it/\left( {re^{i\omega } } \right)} \right|}}\le 1,\; \frac{1}{{\left| {\left( {1 + r\tau e^{i\omega } /it} \right)\left( {1 + it/\left( {re^{i\omega } } \right)} \right)} \right|}} \le 1
\]
and
\[
\frac{1}{{\left| {1 - ise^{i\varphi } /\left( {re^{i\omega } \cos \varphi } \right)} \right|}} \le \sec \left( {\omega  - \varphi } \right),
\]
\[
\frac{1}{{\left| {\left( {1 + ir\tau e^{i\omega } \cos \varphi /\left( {se^{i\varphi } } \right)} \right)\left( {1 - ise^{i\varphi } /\left( {re^{i\omega } \cos \varphi } \right)} \right)} \right|}} \le \sec ^2 \left( {\omega  - \varphi } \right)
\]
for any positive $r$, $s$, $t$ and $\tau$, we deduce the upper bound
\begin{align*}
& \left| {\left( { - 1} \right)^N \frac{{e^{\left( {N + \frac{1}{2}} \right)\left( {\omega  - \theta } \right)i} }}{{\sqrt {2\pi a} i}}\int_0^{ + \infty } {\frac{{\tau ^{N - \frac{1}{2}} e^{ - r\tau \left| {\lambda  + \log \lambda  + 1 + \pi i} \right|} }}{{1 + \tau e^{i\left( {\omega  - \theta } \right)} }}M_K \left( {r\tau e^{i\omega } } \right)d\tau } } \right| \\ & \le \frac{{\sec \left( {\omega  - \varphi } \right)}}{{\cos ^K \varphi }}\frac{1}{{\sqrt {2\pi } r^{K + \frac{1}{2}} }}\left| {\int_0^{ + \infty } {\frac{{\tau ^{N - K - \frac{1}{2}} e^{ - r\tau \left| {\lambda  + \log \lambda  + 1 + \pi i} \right|} }}{{1 + \tau e^{i\left( {\omega  - \theta } \right)} }}d\tau } } \right|\frac{1}{{2\pi }}\int_0^{ + \infty } {s^{K - 1} e^{ - 2\pi s} \left| {\Gamma^\ast \left( {\frac{{ise^{i\varphi } }}{{\cos \varphi }}} \right)} \right|ds} 
\\ & + \frac{{\sec ^2 \left( {\omega  - \varphi } \right)}}{{\cos ^K \varphi }}\frac{{\Gamma \left( {N - K + \frac{1}{2}} \right)}}{{\sqrt {2\pi } \left| {\lambda  + \log \lambda  + 1 + \pi i} \right|^{N - K + \frac{1}{2}} r^{N + 1} }}\frac{1}{{2\pi }}\int_0^{ + \infty } {s^{K - 1} e^{ - 2\pi s} \left| {\Gamma^\ast \left( {\frac{{ise^{i\varphi } }}{{\cos \varphi }}} \right)} \right|ds} 
\\ & + \frac{1}{{\sqrt {2\pi } r^{K + \frac{1}{2}} }}\left| {\int_0^{ + \infty } {\frac{{\tau ^{N - K - \frac{1}{2}} e^{ - r\tau \left| {\lambda  + \log \lambda  + 1 + \pi i} \right|} }}{{1 + \tau e^{i\left( {\omega  - \theta } \right)} }}d\tau } } \right|\frac{1}{{2\pi }}\int_0^{ + \infty } {t^{K - 1} e^{ - 2\pi t} \left| {\Gamma^\ast  \left( { - it} \right)} \right|dt} 
\\ & + \frac{{\Gamma \left( {N - K + \frac{1}{2}} \right)}}{{\sqrt {2\pi } \left| {\lambda  + \log \lambda  + 1 + \pi i} \right|^{N - K + \frac{1}{2}} r^{N + 1} }}\frac{1}{{2\pi }}\int_0^{ + \infty } {t^{K - 1} e^{ - 2\pi t} \left| {\Gamma^\ast \left( { - it} \right)} \right|dt} .
\end{align*}
Employing the inequalities \eqref{eq20}, \eqref{eq21}, the known integral representation of the Riemann Zeta function \cite[25.5.E1]{NIST} and the definition of the terminant function lead to the estimate
\begin{align*}
& \left| {\left( { - 1} \right)^N \frac{{e^{\left( {N + \frac{1}{2}} \right)\left( {\omega  - \theta } \right)i} }}{{\sqrt {2\pi a} i}}\int_0^{ + \infty } {\frac{{\tau ^{N - \frac{1}{2}} e^{ - r\tau \left| {\lambda  + \log \lambda  + 1 + \pi i} \right|} }}{{1 + \tau e^{i\left( {\omega  - \theta } \right)} }}M_K \left( {r\tau e^{i\omega } } \right)d\tau } } \right|
\\ & \le \left( {\frac{{\sec \left( {\omega  - \varphi } \right)}}{{\cos ^K \varphi }} + 1} \right)\left| {e^{a\left( {\lambda  + \log \lambda  + 1 - \pi i} \right)} \sqrt {\frac{{2\pi }}{a}} \widehat T_{N - K + \frac{1}{2}} \left( {a\left( {\lambda  + \log \lambda  + 1 - \pi i} \right)} \right)} \right|\frac{{\zeta \left( K \right)\Gamma \left( K \right)}}{{\left( {2\pi } \right)^{K + 1} \left| a \right|^K }}
\\ & + \left( {\frac{{\sec ^2 \left( {\omega  - \varphi } \right)}}{{\cos ^K \varphi }} + 1} \right)\frac{{\zeta \left( K \right)\Gamma \left( K \right)\Gamma \left( {N - K + \frac{1}{2}} \right)}}{{\left( {2\pi } \right)^{K + \frac{3}{2}} \left| {\lambda  + \log \lambda  + 1 + \pi i} \right|^{N - K + \frac{1}{2}} \left| a \right|^{N + 1} }}.
\end{align*}
Similarly, we have the following upper bound for the other integral in \eqref{eq18}:
\begin{align*}
& \left| {\left( { - 1} \right)^{N+1} \frac{{e^{\left( {N + \frac{1}{2}} \right)\left( {-\omega  - \theta } \right)i} }}{{\sqrt {2\pi a} i}}\int_0^{ + \infty } {\frac{{\tau ^{N - \frac{1}{2}} e^{ - r\tau \left| {\lambda  + \log \lambda  + 1 + \pi i} \right|} }}{{1 + \tau e^{i\left( {-\omega  - \theta } \right)} }}M_K \left( {r\tau e^{-i\omega } } \right)d\tau } } \right|
\\ & \le \left( {\frac{{\sec \left( {\omega  - \varphi } \right)}}{{\cos ^K \varphi }} + 1} \right)\left| {e^{a\left( {\lambda  + \log \lambda  + 1 + \pi i} \right)} \sqrt {\frac{{2\pi }}{a}} \widehat T_{N - K + \frac{1}{2}} \left( {a\left( {\lambda  + \log \lambda  + 1 + \pi i} \right)} \right)} \right|\frac{{\zeta \left( K \right)\Gamma \left( K \right)}}{{\left( {2\pi } \right)^{K + 1} \left| a \right|^K }}
\\ & + \left( {\frac{{\sec ^2 \left( {\omega  - \varphi } \right)}}{{\cos ^K \varphi }} + 1} \right)\frac{{\zeta \left( K \right)\Gamma \left( K \right)\Gamma \left( {N - K + \frac{1}{2}} \right)}}{{\left( {2\pi } \right)^{K + \frac{3}{2}} \left| {\lambda  + \log \lambda  + 1 + \pi i} \right|^{N - K + \frac{1}{2}} \left| a \right|^{N + 1} }}.
\end{align*}
Thus, we conclude that
\begin{align*}
& \left| R_{N,K}\left(a,\lambda\right) \right|
\\ & \le \left( {\frac{{\sec \left( {\omega  - \varphi } \right)}}{{\cos ^K \varphi }} + 1} \right)\left| {e^{a\left( {\lambda  + \log \lambda  + 1 - \pi i} \right)} \sqrt {\frac{{2\pi }}{a}} \widehat T_{N - K + \frac{1}{2}} \left( {a\left( {\lambda  + \log \lambda  + 1 - \pi i} \right)} \right)} \right|\frac{{\zeta \left( K \right)\Gamma \left( K \right)}}{{\left( {2\pi } \right)^{K + 1} \left| a \right|^K }}
\\ & + \left( {\frac{{\sec \left( {\omega  - \varphi } \right)}}{{\cos ^K \varphi }} + 1} \right)\left| {e^{a\left( {\lambda  + \log \lambda  + 1 + \pi i} \right)} \sqrt {\frac{{2\pi }}{a}} \widehat T_{N - K + \frac{1}{2}} \left( {a\left( {\lambda  + \log \lambda  + 1 + \pi i} \right)} \right)} \right|\frac{{\zeta \left( K \right)\Gamma \left( K \right)}}{{\left( {2\pi } \right)^{K + 1} \left| a \right|^K }}
\\ & + \left( {\frac{{\sec ^2 \left( {\omega  - \varphi } \right)}}{{\cos ^K \varphi }} + 1} \right)\frac{{2\zeta \left( K \right)\Gamma \left( K \right)\Gamma \left( {N - K + \frac{1}{2}} \right)}}{{\left( {2\pi } \right)^{K + \frac{3}{2}} \left| {\lambda  + \log \lambda  + 1 + \pi i} \right|^{N - K + \frac{1}{2}} \left| a \right|^{N + 1} }}.
\end{align*}
By continuity, this bound holds in the closed sector $\left| {\arg a} \right| \leq - \omega  + \pi$. As in Appendix \ref{appendixa}, the minimising value $\varphi = \varphi^\ast$ of the factor $ \sec \left( {\omega  - \varphi } \right)\cos ^{ - K} \varphi$ is the unique solution of the equation
\[
\left( {K + 1} \right)\sin \left( {\omega  - 2\varphi^\ast } \right) = \left( {K - 1} \right)\sin \omega ,
\]
that satisfies $ - \frac{\pi}{2} + \omega  < \varphi ^\ast < \frac{\pi }{2} $ if $ \frac{\pi}{2} \le \omega < \pi$, and $0 < \varphi^\ast  < \omega$ if $0  < \omega  < \frac{\pi}{2}$. Assume now that $N = \left| {a\left( {\lambda  + \log \lambda  + 1 + \pi i} \right)} \right| + \rho$ where $\rho$ is bounded. Employing Stirling's formula, we find that
\[
\left( {\frac{{\sec ^2 \left( {\omega  - \varphi } \right)}}{{\cos ^K \varphi }} + 1} \right)\frac{{2\zeta \left( K \right)\Gamma \left( K \right)\Gamma \left( {N - K + \frac{1}{2}} \right)}}{{\left( {2\pi } \right)^{K + \frac{3}{2}} \left| {\lambda  + \log \lambda  + 1 + \pi i} \right|^{N - K + \frac{1}{2}} \left| a \right|^{N + 1} }} = \mathcal{O}_{\lambda ,K,\rho } \left( {\frac{{e^{ - \left| {a\left( {\lambda  + \log \lambda  + 1 + \pi i} \right)} \right|} }}{{\left| {\lambda  + \log \lambda  + 1 + \pi i} \right|^{\frac{1}{2}} \left| a \right|^{K + 1} }}} \right)
\]
as $a\to \infty$. Olver's estimation \eqref{eq22} shows that
\begin{align*}
\left| {e^{a\left( {\lambda  + \log \lambda  + 1 \pm \pi i} \right)} \widehat T_{N - K + \frac{1}{2}} \left( {a\left( {\lambda  + \log \lambda  + 1 \pm \pi i} \right)} \right)} \right| & = \mathcal{O}_{\lambda ,K,\rho } \left( {e^{ - \left| {a\left( {\lambda  + \log \lambda  + 1 \pm \pi i} \right)} \right|} } \right) \\ & = \mathcal{O}_{\lambda ,K,\rho } \left( {e^{ - \left| {a\left( {\lambda  + \log \lambda  + 1 + \pi i} \right)} \right|} } \right)
\end{align*}
for large $a$. Therefore, we obtain that
\begin{equation}\label{eq28}
R_{N,K} \left( {a,\lambda } \right) = \mathcal{O}_{\lambda ,K,\rho } \left( {\frac{{e^{ - \left| {a\left( {\lambda  + \log \lambda  + 1 + \pi i} \right)} \right|} }}{{\left| a \right|^{K + \frac{1}{2}} }}} \right)
\end{equation}
as $a\to \infty$ in the sector $\left| {\arg a} \right| \leq - \omega  + \pi$.

Consider now the sector $ - \omega  + \pi  < \arg a < \omega  + \pi$. When $a$ enters this sector, the pole in the second integral in \eqref{eq18} crosses the integration path. According to the residue theorem, we obtain
\begin{gather}\label{eq27}
\begin{split}
R_{N,K} \left( {a,\lambda } \right) = \; & i\sqrt {\frac{{2\pi }}{a}} e^{a\left( {\lambda  + \log \lambda  + 1 + \pi i} \right)} M_K \left( {ae^{ - \pi i} } \right)
\\ &+ \frac{{\left( { - 1} \right)^N }}{{a^{N + 1} }}\frac{{e^{\left( {N + \frac{1}{2}} \right)\omega i} }}{{\sqrt {2\pi } i}}\int_0^{ + \infty } {\frac{{t^{N - \frac{1}{2}} e^{ - t\left| {\lambda  + \log \lambda  + 1 + \pi i} \right|} }}{{1 + te^{i\omega } /a}}M_K \left( {te^{i\omega } } \right)dt} 
\\ & + \frac{{\left( { - 1} \right)^{N + 1} }}{{a^{N + 1} }}\frac{{e^{ - \left( {N + \frac{1}{2}} \right)\omega i} }}{{\sqrt {2\pi } i}}\int_0^{ + \infty } {\frac{{t^{N - \frac{1}{2}} e^{ - t\left| {\lambda  + \log \lambda  + 1 + \pi i} \right|} }}{{1 + te^{ - i\omega } /a}}M_K \left( {te^{ - i\omega } } \right)dt} 
\end{split}
\end{gather}
for $ - \omega  + \pi  < \arg a < \omega  + \pi$. If $\left| {\arg \left( {ae^{ - \pi i} } \right)} \right| \le \omega  < \pi$, then $M_K \left( {ae^{ - \pi i} } \right) = \mathcal{O}_{\omega ,K} \left( {\left| a \right|^{ - K} } \right) = \mathcal{O}_{\lambda ,K} \left( {\left| a \right|^{ - K} } \right)$ as $a\to \infty$, whence
\begin{equation}\label{eq26}
i\sqrt {\frac{{2\pi }}{a}} e^{a\left( {\lambda  + \log \lambda  + 1 + \pi i} \right)} M_K \left( {ae^{ - \pi i} } \right) = \mathcal{O}_{\lambda ,K} \left( {\frac{{e^{\left| {a\left( {\lambda  + \log \lambda  + 1 + \pi i} \right)} \right|\cos \left( {\arg a + \omega } \right)} }}{{\left| a \right|^{K + \frac{1}{2}} }}} \right)
\end{equation}
as $a\to \infty$ in the sector $ - \omega  + \pi  \leq \arg a \leq \omega  + \pi$. The two integrals can be estimated in the same way as in the case $\left| {\arg a} \right| \leq - \omega  + \pi$, and one finds that they satisfy the order estimate \eqref{eq26} with an implied constant that may also depend on $\rho$. Thus, the final conclusion is that
\[
R_{N,K} \left( {a,\lambda } \right) = \mathcal{O}_{\lambda ,K,\rho } \left( {\frac{{e^{\left| {a\left( {\lambda  + \log \lambda  + 1 + \pi i} \right)} \right|\cos \left( {\arg a + \omega } \right)} }}{{\left| a \right|^{K + \frac{1}{2}} }}} \right)
\]
as $a\to \infty$ in the closed sector $ - \omega  + \pi  \leq \arg a \leq \omega  + \pi$. Similarly, we find that
\[
R_{N,K} \left( {a,\lambda } \right) = \mathcal{O}_{\lambda ,K,\rho } \left( {\frac{{e^{\left| {a\left( {\lambda  + \log \lambda  + 1 + \pi i} \right)} \right|\cos \left( {\arg a - \omega } \right)} }}{{\left| a \right|^{K + \frac{1}{2}} }}} \right)
\]
for large $a$ in the sector $ - \omega  - \pi  \le \arg a \le \omega  - \pi$.

Let $\delta>0$ be a fixed small positive number and suppose that $\delta  < \pi  - \omega$, i.e., $\lambda$ is not very close to $0$. We consider the range $\omega  + \pi  < \arg a \le 2\pi  - \delta  < 2\pi$. To obtain a representation of $R_{N,K} \left( {a,\lambda } \right)$ which is valid in this sector, we rotate the path of integration in the first integral in \eqref{eq27} and apply the residue theorem to find
\begin{align*}
R_{N,K} \left( {a,\lambda } \right) = & - i\sqrt {\frac{{2\pi }}{a}} e^{a\left( {\lambda  + \log \lambda  + 1 - \pi i} \right)} M_K \left( {ae^{ - \pi i} } \right)+ i\sqrt {\frac{{2\pi }}{a}} e^{a\left( {\lambda  + \log \lambda  + 1 + \pi i} \right)} M_K \left( {ae^{ - \pi i} } \right) 
\\ & + \frac{{\left( { - 1} \right)^N }}{{a^{N + 1} }}\frac{{e^{\left( {N + \frac{1}{2}} \right)\omega i} }}{{\sqrt {2\pi } i}}\int_0^{ + \infty } {\frac{{t^{N - \frac{1}{2}} e^{ - t\left| {\lambda  + \log \lambda  + 1 + \pi i} \right|} }}{{1 + te^{i\omega } /a}}M_K \left( {te^{i\omega } } \right)dt} 
\\ & + \frac{{\left( { - 1} \right)^{N + 1} }}{{a^{N + 1} }}\frac{{e^{ - \left( {N + \frac{1}{2}} \right)\omega i} }}{{\sqrt {2\pi } i}}\int_0^{ + \infty } {\frac{{t^{N - \frac{1}{2}} e^{ - t\left| {\lambda  + \log \lambda  + 1 + \pi i} \right|} }}{{1 + te^{ - i\omega } /a}}M_K \left( {te^{ - i\omega } } \right)dt} 
\\ =  & - i\sqrt {\frac{{2\pi }}{a}} e^{a\left( {\lambda  + \log \lambda  + 1 - \pi i} \right)} M_K \left( {ae^{ - \pi i} } \right) + i\sqrt {\frac{{2\pi }}{a}} e^{a\left( {\lambda  + \log \lambda  + 1 + \pi i} \right)} M_K \left( {ae^{ - \pi i} } \right)  + R_{N,K} \left( {ae^{ - 2\pi i} ,\lambda } \right),
\end{align*}
for $\omega  + \pi  < \arg a \le 2\pi  - \delta$. Since $\omega  - \pi  < \arg \left( {ae^{ - 2\pi i} } \right) \le  - \delta$, $R_{N,K} \left( {ae^{ - 2\pi i} ,\lambda } \right)$ satisfies the order estimate given in the right-hand side of \eqref{eq28}. If $\omega  < \arg \left( {ae^{ - \pi i} } \right) \le \pi  - \delta$, then $M_K \left( {ae^{ - \pi i} } \right) = \mathcal{O}_{\omega ,K,\delta} \left( {\left| a \right|^{ - K} } \right) = \mathcal{O}_{\lambda ,K,\delta} \left( {\left| a \right|^{ - K} } \right)$ as $a\to \infty$, whence
\begin{multline*}
 - i\sqrt {\frac{{2\pi }}{a}} e^{a\left( {\lambda  + \log \lambda  + 1 - \pi i} \right)} M_K \left( {ae^{ - \pi i} } \right) + i\sqrt {\frac{{2\pi }}{a}} e^{a\left( {\lambda  + \log \lambda  + 1 + \pi i} \right)} M_K \left( {ae^{ - \pi i} } \right)
\\ = \mathcal{O}_{\lambda ,K,\delta } \left( {\frac{{e^{\left| {a\left( {\lambda  + \log \lambda  + 1 + \pi i} \right)} \right|\cos \left( {\arg a - \omega } \right)}  + e^{\left| {a\left( {\lambda  + \log \lambda  + 1 + \pi i} \right)} \right|\cos \left( {\arg a + \omega } \right)} }}{{\left| a \right|^{K + \frac{1}{2}} }}} \right)
\end{multline*}
as $a\to \infty$ in the sector $\omega  + \pi  < \arg a \le 2\pi  - \delta$. It is elementary to show that in this range $\cos \left( {\arg a - \omega } \right) < \cos \left( {\arg a + \omega } \right)$, and therefore, we conclude that
\[
R_{N,K} \left( {a,\lambda } \right) = \mathcal{O}_{\lambda ,K,\rho,\delta } \left( {\frac{{e^{\left| {a\left( {\lambda  + \log \lambda  + 1 + \pi i} \right)} \right|\cos \left( {\arg a + \omega } \right)} }}{{\left| a \right|^{K + \frac{1}{2}} }}} \right)
\]
as $a\to \infty$ in the sector $\omega  + \pi  < \arg a \le 2\pi  - \delta$. Similarly, we find that
\[
R_{N,K} \left( {a,\lambda } \right) = \mathcal{O}_{\lambda ,K,\rho,\delta } \left( {\frac{{e^{\left| {a\left( {\lambda  + \log \lambda  + 1 + \pi i} \right)} \right|\cos \left( {\arg a - \omega } \right)} }}{{\left| a \right|^{K + \frac{1}{2}} }}} \right)
\]
for large $a$ in the sector $-2\pi + \delta \le \arg a < -\pi - \omega$.

Consider finally the cases $K = 0$ and $K = 1$. We can write
\begin{align*}
R_{N,0} \left( {a,\lambda } \right) = & - ie^{a\left( {\lambda  + \log \lambda  + 1 - \pi i} \right)} \sqrt {\frac{{2\pi }}{a}} \sum\limits_{k = 0}^1 {\frac{{\gamma _k }}{{a^k }}\widehat T_{N - k + \frac{1}{2}} \left( {a\left( {\lambda  + \log \lambda  + 1 - \pi i} \right)} \right)} 
\\ & + ie^{a\left( {\lambda  + \log \lambda  + 1 + \pi i} \right)} \sqrt {\frac{{2\pi }}{a}} \sum\limits_{k = 0}^1 {\frac{{\gamma _k }}{{a^k }}\widehat T_{N - k + \frac{1}{2}} \left( {a\left( {\lambda  + \log \lambda  + 1 + \pi i} \right)} \right)} 
\\ & + R_{N,2} \left( {a,\lambda } \right)
\end{align*}
and
\begin{align*}
R_{N,1} \left( {a,\lambda } \right) = & - ie^{a\left( {\lambda  + \log \lambda  + 1 - \pi i} \right)} \sqrt {\frac{{2\pi }}{a}} \frac{{\gamma _1 }}{a}\widehat T_{N - \frac{1}{2}} \left( {a\left( {\lambda  + \log \lambda  + 1 - \pi i} \right)} \right)
\\ & + ie^{a\left( {\lambda  + \log \lambda  + 1 + \pi i} \right)} \sqrt {\frac{{2\pi }}{a}} \frac{{\gamma _1 }}{a}\widehat T_{N - \frac{1}{2}} \left( {a\left( {\lambda  + \log \lambda  + 1 + \pi i} \right)} \right)
\\ & + R_{N,2} \left( {a,\lambda } \right).
\end{align*}
Employing the previously obtained bounds for $R_{N,2} \left( {a,\lambda } \right)$ and Olver's estimation \eqref{eq22} together with the connection formula for the terminant function \cite[p. 260]{Paris2}, shows that $R_{N,0} \left( {a,\lambda } \right)$ and $R_{N,1} \left( {a,\lambda } \right)$ indeed satisfy the order estimates prescribed in Theorem \ref{thm2}.

\appendix

\section{An auxiliary estimate}\label{appendixa} In this appendix, we derive a bound for the remainder term $M_N \left( z \right)$ in the asymptotic series of the scaled gamma function, when $0<\arg z<\pi$ and $N\geq 2$. Let $0<\varphi<\frac{\pi}{2}$ be an arbitrary acute angle that satisfies $\varphi <\arg z<\frac{\pi}{2}+\varphi$. We rotate the path of integration in the first integral in \eqref{eq14} and perform the change of variable $t = \frac{{se^{i\varphi } }}{\cos \varphi}$, to obtain the representation
\begin{equation}\label{eq19}
M_N \left( z \right) = \frac{1}{{2\pi i}}\frac{{i^N }}{{z^N }}\left( {\frac{{e^{i\varphi } }}{{\cos \varphi }}} \right)^N \int_0^{ + \infty } {\frac{{s^{N - 1} e^{ - 2\pi \frac{{se^{i\varphi } }}{{\cos \varphi }}} \Gamma^\ast \left( {\frac{{ise^{i\varphi } }}{{\cos \varphi }}} \right)}}{{1 - ise^{i\varphi } /z\cos \varphi }}ds}  - \frac{1}{{2\pi i}}\frac{{\left( { - i} \right)^N }}{{z^N }}\int_0^{ + \infty } {\frac{{t^{N - 1} e^{ - 2\pi t} \Gamma^\ast \left( { - it} \right)}}{{1 + it/z}}dt} ,
\end{equation}
which is valid when $\varphi <\arg z<\frac{\pi}{2}+\varphi$, using analytic continuation. Simple estimation and the inequality \eqref{eq40} yield
\begin{equation}\label{eq15}
\left| {M_N \left( z \right)} \right| \le \frac{1}{{2\pi }}\frac{1}{{\left| z \right|^N }}\frac{{\sec \left( {\theta  - \varphi } \right)}}{{\cos ^N \varphi }}\int_0^{ + \infty } {s^{N - 1} e^{ - 2\pi s} \left| {\Gamma^\ast \left( {\frac{{ise^{i\varphi } }}{{\cos \varphi }}} \right)} \right|ds}  + \frac{1}{{2\pi }}\frac{1}{{\left| z \right|^N }}\int_0^{ + \infty } {t^{N - 1} e^{ - 2\pi t} \left| {\Gamma^\ast \left( { - it} \right)} \right|dt} ,
\end{equation}
with the notation $\theta = \arg z$. It was proved in \cite{Nemes} that for any $s>0$ and $0<\varphi<\frac{\pi}{2}$, it holds that
\begin{equation}\label{eq20}
\left| {\Gamma^\ast \left( {\frac{{ise^{i\varphi } }}{{\cos \varphi }}} \right)} \right| \le \frac{1}{{1 - e^{ - 2\pi s} }}.
\end{equation}
We also have
\begin{equation}\label{eq21}
\left| {\Gamma^\ast \left( { - it} \right)} \right| = \frac{1}{{\sqrt {1 - e^{ - 2\pi t} } }} \le \frac{1}{{1 - e^{ - 2\pi t} }}.
\end{equation}
Substituting these estimates into \eqref{eq15} gives the bound
\begin{equation}\label{eq16}
\left| {M_N \left( z \right)} \right| \le \left( {\frac{{\sec \left( {\theta  - \varphi } \right)}}{{\cos ^N \varphi }} + 1} \right)\frac{{\zeta \left( N \right)\Gamma \left( N \right)}}{{\left( {2\pi } \right)^{N + 1} \left| z \right|^N }},
\end{equation}
for $N\geq 2$ and $\varphi <\arg z<\frac{\pi}{2}+\varphi$. Here we have made use of the known integral representation of the Riemann Zeta function \cite[25.5.E1]{NIST}. The minimisation of the factor $ \sec \left( {\theta  - \varphi } \right)\cos ^{ - N} \varphi$ as a function of $\varphi$ can be done using a lemma of Meijer \cite[pp. 953--954]{Meijer}. In our case, Meijer's lemma gives that the minimising value $\varphi = \varphi^\ast$ in \eqref{eq16}, is the unique solution of the equation
\[
\left( {N + 1} \right)\sin \left( {\theta  - 2\varphi^\ast } \right) = \left( {N - 1} \right)\sin \theta ,
\]
that satisfies $ - \frac{\pi}{2} + \theta  < \varphi ^\ast < \frac{\pi }{2} $ if $ \frac{\pi}{2} \le \theta < \pi$, and $0 < \varphi^\ast  < \theta$ if $0  < \theta  < \frac{\pi}{2}$. With this choice of $\varphi$, \eqref{eq16} provides a bound for $M_N \left( z \right)$ in the range $0<\arg z<\pi$.

\section*{Acknowledgment} The research of the author was supported by the Central European University Foundation, Budapest (CEUBPF). The author would like to thank Adri B. Olde Daalhuis for his useful comments and suggestions on the manuscript during the author's visit at the University of Edinburgh.

\end{document}